\documentclass{amsart}
\usepackage{amssymb}
\usepackage{amscd}
\usepackage[all]{xy}
\newcommand\datver[1]{\def\datverp%
 {\par\boxed{\boxed{\text{#1; Run: \today}}}}}
\datver{Revised: 16-8-2005}
\numberwithin{equation}{section}
\newtheorem{theorem}{Theorem}[section]
\newtheorem{proposition}[theorem]{Proposition}
\newtheorem{corollary}[theorem]{Corollary}
\newtheorem{definition}[theorem]{Definition}
\newtheorem{lemma}[theorem]{Lemma}
\newtheorem{remark}[theorem]{Remark}
\newtheorem{notation}[theorem]{Notation}

\renewcommand\Re{\operatorname{Re}}

\newcommand\bbC{\mathbb C}

\newcommand\bbN{\mathbb N}

\newcommand\bbR{\mathbb R}
\newcommand\bbS{\mathbb S}

\newcommand\bbZ{\mathbb Z}
\newcommand\cb{\overline{B}^{k}_{1}}
\newcommand\CI{{\mathcal{C}}^{\infty}}
\newcommand\CId{\dot{\mathcal{C}}^{\infty}}
\newcommand\coker{\operatorname{coker}}
\newcommand\E{\mathcal{E}}
\newcommand\F{\mathcal{F}}
\newcommand\K{K^{0}_{c}}
\newcommand\Gv{\mathcal{G}}
\newcommand\Hv{\mathcal{H}}
\newcommand\End{\operatorname{End}}
\newcommand\fc{\Psi_{\Phi}}
\newcommand\fci{\Psi_{\Id}}
\newcommand\fcde{{}^{\Phi}\Omega}
\newcommand\fcd{\operatorname{Diff}_{\Phi}}
\newcommand\fcn{{}^{\Phi}N}
\newcommand\fcis{{}^{\Id}N}
\newcommand\fcs{{}^{\Phi}S}
\newcommand\fct{{}^{\Phi}T}
\newcommand\fcv{\mathcal{V}_{\Phi}}
\newcommand\ff{\operatorname{ff}}
\newcommand\Frc{\mathcal{F}^{-\infty}_{\Phi}}
\newcommand\Frco{\mathcal{F}^{-\infty}_{\Phi,0}}
\newcommand\Fri{\mathcal{F}^{-\infty}_{\Id}}
\newcommand\fs{\Psi_{\Phi s}}
\newcommand\G{G^{-\infty}}
\newcommand\Gc{G^{-\infty}_{\Phi}}
\newcommand\Gs{G^{-\infty}_{\Phi s}}
\newcommand\Gso{G^{-\infty}_{\Phi s, 0}}
\newcommand\Gd{\dot{G}^{-\infty}}
\newcommand\GL{\operatorname{GL}}
\newcommand\Hom{\operatorname{Hom}}
\newcommand\Hc{\operatorname{H}_{\Phi}}
\newcommand\ind{\operatorname{ind}}
\newcommand\Id{\operatorname{Id}}
\newcommand\Ld{\operatorname{L}^{2}}
\newcommand\no{N_{\Phi}}
\newcommand\noi{N_{\Id}}
\newcommand\nos{N_{\Phi *}}
\newcommand\pa{\partial}
\newcommand\pt{\operatorname{pt}}
\newcommand\rank{\operatorname{rank}}
\newcommand\supp{\operatorname{supp}}
\newcommand\Tr{\operatorname{Tr}}
\newcommand\tr{\operatorname{tr}}
\newcommand\rTr{\operatorname{Tr_{R}}}
\newcommand\bTr{\overline{\operatorname{Tr}}}
\newcommand\Vy{\mathcal{V}_{Y}}
\newcommand\W{\mathcal{W}_{Y}}

\begin{document}

\title[Bott Periodicity for Fibred Cusp Operators]
{Bott Periodicity for Fibred Cusp Operators}
\author{Fr\'ed\'eric Rochon}
\thanks{The author was partially supported
by Le Fonds Qu\'{e}b\'{e}cois de la recherche sur la nature et les
technologies and the Natural Sciences and Engineering Research Council of 
Canada.} 
\address{Department of Mathematics, Massachusetts Institute of Technology}
\email{rochon@math.mit.edu}
\begin{abstract}
In the framework of fibred cusp operators on a manifold $X$ associated to 
a boundary fibration $\Phi: \pa X\to Y$, the homotopy groups of the 
space $\Gc(X;E)$ of invertible smoothing perturbations of the identity 
are computed in terms of the $K$-theory of $T^{*}Y$.  It is shown that there
is a periodicity, namely the odd and the even homotopy groups are isomorphic
among themselves.  To obtain this result, one of the important 
steps is the description of the index of a Fredholm smoothing 
perturbation of the identity in terms of an associated $K$-class
in $\K(T^{*}Y)$. 
\end{abstract}
\maketitle


\section*{Introduction}\label{int.0}

For standard pseudodifferential operators on a closed manifold $X$ acting
on some complex vector bundle $E$, Bott periodicity arises by considering
the group
\[
   \G(X;E)=\{ \Id +Q \; | \; Q\in \Psi^{-\infty}(X;E), \quad \Id+Q \;
\mbox{is invertible} \}
\]
of invertible smoothing perturbations of the identity.  This becomes a 
topological group by taking the $\CI$-topology induced by the 
identification of smoothing operators with their Schwartz kernels, which
are smooth sections of some bundle over $X\times X$.  If $\Delta^{E}$ is any
Laplacian acting on sections of $E$, and if $\{f_{i}\}_{i\in\bbN}$ is a basis
of $\Ld(X;E)$ coming from a sequence of orthonormal eigensections of 
$\Delta^{E}$ with increasing eigenvalues, then there is an isomorphism of 
topological groups
\[
\begin{array}{rrcl}
     f_{\Delta^{E}}: & \G(X;E) &\to & \Gv^{-\infty} \\
                     & \Id +Q & \mapsto & \delta_{ij}+ \langle f_{i},Qf_{j}
\rangle,
\end{array}
\]
where $\Gv^{-\infty}$ is the group of invertible semi-infinite matrices 
$\delta_{ij}+Q_{ij}$ such that
\[
              \|Q\|_{k}= \sum_{i,j}(i+j)^{k}|Q_{ij}|< \infty, \quad \forall
k\in \bbN_{0},
\]
the topology of $\Gv^{-\infty}$ being the one induced by the norms 
$\| \cdot \|_{k}$, $k\in\bbN_{0}$. 

This isomorphism indicates that the topology of $\G(X;E)$ does not depend
at all on the geometry of $X$ and $E$.  Since the direct limit
\[
        \GL(\infty,\bbC)= \lim_{k\to \infty} \GL(k,\bbC)
\] 
is a weak deformation retract (see definition~\ref{hg.5} below) of 
$\Gv^{-\infty}$, they share the same homotopy groups, which is to 
say 
\[
         \pi_{k}(\G(X;E))\cong \pi_{k}(\Gv^{-\infty})\cong
\left\{ \begin{array}{cl}
  \{0\} & \quad k\;\mbox{even}, \\
    \bbZ & \quad k\;\mbox{odd}.
\end{array} \right.
\]
This periodicity in the homotopy groups is an instance of Bott periodicity 
as originally described by Bott in \cite{Bott}.

In this paper, we will describe how Bott periodicity arises when one considers
instead fibred cusp operators $\fc^{*}(X;E)$ on a compact manifold with 
boundary $X$ acting on some complex vector bundle $E$.  These operators were
introduced by Mazzeo and Melrose in \cite{mazzeo-melrose4}.  The definition
involves a defining function for the boundary $\pa X$ and a fibration 
$\Phi:\pa X\to Y$ of the boundary.  Again, one considers the group
\[
    \Gc(X;E)=\{ \Id +Q\; |\; Q\in \fc^{-\infty}(X;E),\quad \Id +Q \;
\mbox{is invertible}\}
\]
of invertible smoothing perturbations of the identity.  As before, 
it has a $\CI$-topology induced by the identification 
of smoothing operators with their Schwartz kernels.  \textbf{Our
main result, stated in theorem~\ref{hg.31}}, is to describe the homotopy groups
of $\Gc(X;E)$ in terms of the $K$-theory of $T^{*}Y$, namely, all the even
homotopy groups are shown to be isomorphic to the kernel of the 
topological index map 
\[
\ind_{t}:\K(T^{*}Y)\to \bbZ,
\]
while all the odd
homotopy groups are shown to be isomorphic to $\widetilde{K}^{-1}(Y^{T^{*}Y})$,
where $Y^{T^{*}Y}$ is the Thom space associated to the vector bundle $T^{*}Y$.
This periodicity of the homotopy groups is what we interpret as Bott
periodicity for fibred cusp operators.  Strictly speaking, this result is only
true when the fibres of the fibration $\Phi:\pa X\to Y$ are of dimension at 
least one, but in the particular case where $\Phi:\pa X\to Y$ is a finite
covering, which includes the case of scattering operators, the result is 
still true provided one allows some stabilization (see the discussion at the
beginning of section~\ref{hg.0}).

This is in a certain sense a generalization of proposition $3.6$ in 
\cite{fipomb}, where
it was shown that, in the particular case where $\Phi:\pa X\to \pt$ is a 
trivial fibration (the case of cusp operators), all the homotopy groups
of $\Gc(X;E)$ are trivial.  This weak contractibility was used in 
\cite{fipomb} to derive a relative index theorem for families of elliptic
cusp pseudodifferential operators.  However, this relative index theorem 
does not
seem to generalize in a simple way to fibred cusp operators with a non-trivial
fibration $\Phi:\pa X\to Y$, but we still hope that the knowledge of the 
homotopy groups of $\Gc(X;E)$ will turn out to be useful in the understanding
of the index of general Fredholm fibred cusp operators.

An important feature of fibred cusp operators is that smoothing operators 
are not necessarily compact.  As a consequence, smoothing perturbations of 
the identity are not necessarily Fredholm, and when they are, they do not
necessarily have a vanishing index.  Our computation of the homotopy groups
of $\Gc(X;E)$ relies on a careful study of the space 
\[
   \Frc(X;E)=\{ \Id+Q \; | \; Q\in \fc^{-\infty}(X;E)\quad \Id +Q\;
\mbox{is Fredholm}\}
\]   
of Fredholm smoothing perturbations of the identity.  
\textbf{The second important
result of this paper, stated in theorem~\ref{im.9}}, is that the index
of a Fredholm operator $(\Id+Q)\in\Frc(X;E)$ can be described in terms
of an associated $K$-class $\kappa(\Id+Q)\in\K(T^{*}Y)$, namely
\[
         \ind(\Id+Q)=\ind_{t}(\kappa(\Id+Q)),
\]
where $\ind_{t}$ is the topological index introduced by Atiyah and Singer
in \cite{Atiyah-Singer1}.  When $\Phi:\pa X\to Y$ is not a finite covering, an
important corollary is that the index map 
\[
        \ind: \Frc(X;E)\to \bbZ
\]
is surjective.  However, let us emphasize that the topology of $\Frc(X;E)$ 
is in general significantly different from the topology of the space
$\mathcal{F}(X;E)$ of all Fredholm operators acting on $\Ld(X;E)$.  Indeed,
$\mathcal{F}(X;E)$, which is a classifying space for even $K$-theory, has 
trivial odd homotopy groups and even homotopy groups isomorphic to $\bbZ$, 
while $\Frc(X;E)$, as described in proposition~\ref{hg.3} below, has even
homotopy groups isomorphic to $\K(T^{*}Y)$ and odd homotopy groups
isomorphic to $\widetilde{K}^{-1}(Y^{T^{*}Y})$.  Nevertheless, in the 
particular case of cusp operators ($Y=\pt$), this shows that
$\Frc(X;E)\subset \mathcal{F}(X;E) $ is quite big and is also a classifying
space for even $K$-theory.

The paper is organized as follows.  We first briefly review the definition and
the main properties of fibred cusp operators in section~\ref{fco.0}.  We then
recall in section~\ref{tdf.0} the notion of a regularized trace on smoothing
fibred cusp operators and its associated trace-defect formula.  This is 
used in section~\ref{itdf.0} to get an analytic formula for the index
of operators in $\Frc(X;E)$.  Section~\ref{ss.0} is about spectral sections,
which play a crucial r\^{o}le in section~\ref{ckc.0}, where a $K$-class 
is associated to any Fredholm operator in
$\Frc(X;E)$.  This $K$-class is in turn used in section~\ref{im.0} to get
a topological index formula.  This will allow us, in section~\ref{hg.0}, to
present the main result of this paper, that is, the periodicity of the 
homotopy groups of $\Gc(X;E)$.  Finally, in section~\ref{vhg.0}, we discuss
the relation of our result with \cite{fipomb} and show that, once the 
weak contractibility is known, for instance in the case of cusp operators, it
is not much harder, adapting a proof of Kuiper in \cite{Kuiper}, to deduce the
actual contractibility, this last improvement being more esthetic than 
useful.

\section*{Acknowledgments}

The content of this paper is included in the 
Ph.D. thesis of the author.  The author is very grateful to his Ph.D.
advisor Richard B. Melrose, who was very generous of his time and 
ideas.  The author would like to thank Andr\'{a}s Vasy for useful comments
on the manuscript.  
The author would like also to thank Benoit Charbonneau, K\'{a}ri 
Ragnarsson
and Nata\v{s}a \v{S}e\v{s}um for helpful discussions.

\section{Fibred Cusp Operators}\label{fco.0}

Since fibred cusp operators are the central object of study in this paper, let 
us recall briefly their definition and their main properties.  A detailed 
description can be found in \cite{mazzeo-melrose4},
where they were originally introduced.  One can also look at section $2$ of 
\cite{Lauter-Moroianu}.  For a more general discussion in terms of 
Lie algebroids, we refer to \cite{Ammann-Lauter-Nistor} and \cite{Nest-Tsygan}.

Let $X$ be a compact connected manifold with non-empty boundary $\pa X$.  
Assume that the
boundary $\pa X$ has the structure of a (locally trivial) fibration
\begin{equation}
\xymatrix{Z\ar@{-}\ar[r]&\pa X\ar[d]^{\Phi}\\&Y}
\label{fco.1}\end{equation}
where $Y$ and $Z$ are closed manifolds, $Y$ being the base of the fibration and
$Z$ being a typical fibre.  We will not assume that the boundary $\pa X$ is 
connected, but if $Y$ is disconnected, we will assume that $Z$ is the same
over all connected components of $Y$.  For the disconnected case, we will 
basically use the 
generalization of \cite{fipomb} rather than the one presented in 
\cite{mazzeo-melrose4} (see remark~\ref{fco.9} below).  Let $x\in\CI(X)$ be a 
defining function for the boundary $\pa X$, that is, $x$ is non-negative, 
$\pa X$ is the zero set of $x$ and $dx$ does not vanish anywhere on $\pa X$.

The defining function $x$ gives rise to an explicit trivialization of the 
conormal bundle of $\pa X$.  A fibration \eqref{fco.1} together with 
an explicit trivialization of the conormal bundle of $\pa X$ defines a 
fibred cusp structure.  To a given fibred cusp structure, we will associate
 an algebra of fibred cusp operators denoted $\fc^{*}(X)$.  Changing the 
trivialization of the conormal bundle of $\pa X$ will not change $\fc^{*}(X)$
drastically, since we will get an isomorphic algebra, but not in a 
canonical way.  When we talk about an algebra of fibred cusp operators, 
we
therefore tacitly assume that a defining function has been chosen.  On the
other hand, changing the fibration \eqref{fco.1} has important consequences
on the algebra $\fc^{*}(X)$.  In particular, there are two extreme cases:
the case where $Z=\{\pt\}$ and $\Phi=\Id$, which leads to the algebra of
\textbf{scattering operators}, and the case where $Y=\{\pt\}$, which leads to 
the algebra of \textbf{cusp operators}. 

\begin{definition}
A \textbf{fibred cusp vector field} $V\in\CI(X, TX)$ is a vector field which, 
at the boundary $\pa X$, is tangent to the fibres of the fibration $\Phi$, 
and such that $Vx\in x^{2}\CI(X)$, where $x$ is the defining function of the 
fibred cusp structure.  We denote by $\fcv(X)$ the Lie algebra of fibred cusp 
vector fields.  
\label{fco.2}\end{definition} 
  
Let $(x,y,z)$ be coordinates in a neighborhood of $p\in\pa X\subset X$, where 
$x$ is the defining function for $\pa X$ and $y$ and $z$ are respectively 
local coordinates on $Y$ and $Z$ (we assume that the fibration $\Phi$ is 
trivial on the neighborhood of $p$ we consider).  Then any fibred cusp
vector field $V\in \fcv(X)$ is locally of the form
\begin{equation}
V= a x^{2}\frac{\pa}{\pa x}+ \sum_{i=1}^{l} b_{i}x\frac{\pa}{\pa y^{i}}
            +\sum_{i=1}^{k} c_{i}\frac{\pa}{\pa z^{i}},
\label{fco.3}\end{equation}  
where $l= \dim Y$, $k= \dim Z$ and $a,b_{i}, c_{i}$ are smooth functions.

\begin{definition}
The space of \textbf{fibred cusp differential operators} of order $m$, denoted 
$\fcd^{m}(X)$, is the space of operators on $\CI(X)$ generated by $\CI(X)$ and
products of up to $m$ elements of $\fcv(X)$.
\label{fco.4}\end{definition}

In the local coordinates $(x,y,z)$, a fibred cusp differential operator
$P\in \fcd^{m}(X)$ may be written as
\begin{equation}
P= \sum_{|\alpha| + |\beta| +q\le m} p_{\alpha,\beta,q} 
(x^{2}\frac{\pa}{\pa x})^{q}(x\frac{\pa}{\pa y})^{\alpha}
(\frac{\pa}{\pa z})^{\beta}
\label{fco.5}\end{equation}
where the $p_{\alpha,\beta,q}$ are smooth functions.

Intuitively, fibred cusp pseudodifferential operators are a generalization
of definition~\ref{fco.4} by allowing $P$ to be not only a polynomial 
in $x^{2}\frac{\pa}{\pa x},x\frac{\pa}{\pa y},\frac{\pa}{\pa z}$, but also
a more general ``function'' in $x^{2}\frac{\pa}{\pa x},x\frac{\pa}{\pa y},
\frac{\pa}{\pa z}$.  The rigorous definition is however better understood by
describing the Schwartz kernels of fibred cusp pseudodifferential operators.

Consider the cartesian product $X^{2}$ of the manifold $X$ with itself.  This 
is a manifold with corner, the corner being $\pa X\times \pa X\subset X^{2}$.
Schwartz kernels of fibred cusp pseudodifferential operators  will be 
distributions on $X^{2}$.  But all the subtleties of their behavior turn
out to happen in the corner of $X^{2}$.  To have a better picture of 
what is going on, one therefore blows up the corner in the sense of Melrose,
\cite{MelroseMWC}, \cite{MelroseAPS}.  In fact, we need to blow up twice. 

We first blow up the corner to get the space  $X^{2}_{b}=[X^{2};\pa X\times 
\pa X]$.  In \cite{MelroseAPS}, this space is used to describe
the Schwartz kernels of $b$-pseudodifferential operators, where ``$b$'' stands
for boundary.  For fibred cusp operators, we need to perform a second 
blow-up which depends on the fibration $\Phi$ we consider.

Beside the two old boundaries coming from $X^{2}$, $X^{2}_{b}$ has a new 
boundary called the \textbf{front face}.  If $\beta_{b}:X^{2}_{b}\to X^{2}$
is the blow-down map, then this front face, denoted $\ff_{b}$, is given by
\begin{equation}
\ff_{b}=\beta_{b}^{-1}(\pa X\times \pa X)\cong (\pa X\times \pa X)\times 
[-1,1]_{s},
\label{fco.7}\end{equation}  
where $s=\frac{x-x'}{x+x'}$ and $x,x'$ are the same boundary defining functions
on the left and right factors of $X^{2}$.  The function $s$ is well-defined 
on $X^{2}_{b}$.  What we blow up next is the submanifold 
$B_{\Phi}\subset\ff_{b}\cong (\pa X\times \pa X)\times [-1,1]_{s}$ given by
\begin{equation}
B_{\Phi}=\{ (h,h',0)\in \ff_{b}\; | \; \Phi(h)=\Phi(h')\}.
\label{fco.8}\end{equation} 
\begin{remark}
When $\pa X$ is disconnected, we use the same definition for $B_{\Phi}$ and
$\ff_{b}$. This is different from \cite{mazzeo-melrose4}, where they consider 
smaller $B_{\Phi}$ and $\ff_{b}$ in the disconnected case.  Using
the same formal definition as in the connected case makes the generalization to
the disconnected case straightforward. 
\label{fco.9}\end{remark}
So we consider the space $X^{2}_{\Phi}=[X^{2}_{b};B_{\Phi}]$.  If 
$\beta_{\Phi}:X^{2}_{\Phi}\to X^{2}_{b}$ is the blow-down map, then the new
boundary appearing on $X^{2}_{\Phi}$ is given by 
\begin{equation}
\ff_{\Phi}=\beta_{\Phi}^{-1}(B_{\Phi}).
\label{fco.10}\end{equation} 

If $(y,z)$ and $(y',z')$ are local coordinates on the left and right factors
of $(\pa X)^{2}$, then
\begin{equation}
S=\frac{x-x'}{x^{2}}, Y=\frac{y-y'}{x}, z-z', x, y, z
\label{fco.11}\end{equation}
are local coordinates on $X^{2}_{\Phi}$, and in these coordinates, $\ff_{\Phi}$
occurs where $x=0$.  Under the blow-down maps $\beta_{b}$ and $\beta_{\Phi}$,
the diagonal $\Delta\subset \pa X\times \pa X$ lifts to $\Delta_{\Phi}\subset
X^{2}_{\Phi}$.  In the coordinates \eqref{fco.11}, it occurs where
$S=Y=z-z'=0$.  

\begin{definition}
The \textbf{fibred cusp density bundle} $\fcde$ is the space of densities on 
$X$ which are locally of the form
$\frac{\alpha}{x^{l+2}}dx dy dz$ near the boundary $\pa X$, where $\alpha$ is a
smooth function and $l=\dim Y$.  Let $\fcde_{R}$ be the lift of $\fcde$ to
$X^{2}$ from the right factor.  This gives the 
\textbf{right fibred cusp density
bundle} $\fcde_{R}'=\beta_{\Phi}^{*}\beta_{b}^{*}(\fcde_{R})$ on 
$X^{2}_{\Phi}$.
\label{fco.21}\end{definition}

\begin{definition}
For any $m\in \bbR$, the space of $\Phi$-pseudodifferential operators of 
order $m$ is given by
\[
       \fc^{m}(X)=\{ K\in I^{m}(X^{2}_{\Phi},\Delta_{\Phi}; \fcde_{R}')\; | \;
              K\equiv 0 \hspace{0.2cm} \mbox{at} \hspace{0.2cm} 
             \pa X^{2}_{\Phi} \setminus \ff_{\Phi}\}
\]
where $I^{m}(X^{2}_{\Phi},\Delta_{\Phi}; \fcde_{R}')$ denotes the space of 
conormal distributions (to $\Delta_{\Phi}$) of order $m$ (we refer to 
\cite{mazzeo-melrose4} and \cite{hormander3} for more details) and $\equiv$
denotes equality in Taylor series.
\label{fco.12}\end{definition}

In local coordinates near the boundary, the action (near the diagonal)
of fibred cusp pseudodifferential operators can be described as follows.

\begin{proposition}
If $\chi\in\CI(X)$ has support in a coordinate patch $(\mathcal{U},(x,y,z))$ 
with $\mathcal{U}\cap\pa X\ne \emptyset $, where $x$ is the defining function 
of 
$\pa X$ and $y$ and $z$ are coordinates on $Y$ and $Z$, then the action of
$P\in \fc^{m}(X)$  on $u\in \CI_{c}(\mathcal{U})$ can be written  
\[
   (\chi Pu)(x,y,z)=\int P_{\chi}(x,y,z,S,Y,z-z')\tilde{u}(x(1-xS),y-xY,z')
     dS dY dz'
\]
where $S=\frac{x-x'}{x^{2}}$, $Y=\frac{y-y'}{x}$, $\tilde{u}$ is the coordinate
representation of $u$ and $P_{\chi}$ is the restriction to $\mathcal{U}
\times \bbR^{n}$ of a distribution on $\bbR^{n}\times\bbR^{n-k}\times
\bbR^{k}$ which has compact support in the first and third variables, is
conormal to $\{S=0, Y=0\}\times\{z=z'\}$ and is rapidly decreasing with 
all derivatives as $|(S,Y)|\to\infty$.  
\label{fco.22}\end{proposition}

If $E$ and $F$ are complex vector bundles on $X$, one can easily extend this
definition to $\Phi$-pseudodifferential operators $\fc^{m}(X;E,F)$ of order
$m$ which act from $\CI(X;E)$ to $\CI(X;F)$.

\begin{proposition}\textbf{(Composition)} For $m,n\in \bbR$ and $E,F,G$ complex
vector bundles on $X$, we have that
\[
         \fc^{m}(X;F,G)\circ\fc^{n}(X;E,F)\subset \fc^{m+n}(X;E,G).
\]
\label{fco.13}\end{proposition} 

As in the case of usual pseudodifferential operators, there is a notion
of ellipticity.

\begin{definition}
The \textbf{fibred cusp tangent bundle} $\fct X$ is the smooth vector bundle
on $X$ such that $\fcv(X)=\CI(X; \fct X)$. 
\label{fco.14}\end{definition}
Notice that $\fct X$ is isomorphic to $TX$, although not in a natural way.  Let
$\fcs^{*}X=(\fct X\setminus 0)/\bbR^{+}$ denote the fibred cusp sphere bundle.
Let $\mathcal{R}^{m}$ be the trivial complex line bundle on $\fcs^{*}X$ with
sections given by functions over $\fct X\setminus 0$ which are positively 
homogeneous of degree $m$.

\begin{proposition}
For all $m\in \bbZ$, there exists a symbol map 
\[
           \sigma_{m}:\fc^{m}(X;E,F)\to \CI(\fcs^{*}X;\hom(E,F)\otimes
                                            \mathcal{R}^{m})
\]
such that we have the following short exact sequence
\[
    \fc^{m-1}(X;E,F)\longrightarrow \fc^{m}(X;E,F)
            \overset{\sigma_{m}}\longrightarrow \CI(\fcs^{*}X;\hom(E,F)\otimes
                                            \mathcal{R}^{m}).     
\]
\label{fco.23}\end{proposition}

\begin{definition}
If $P\in \fc^{m}(X;E,F)$, we say that $\sigma_{m}(P)$ is the \textbf{symbol} 
of $P$.  Moreover, if the symbol $\sigma_{m}(P)$ is invertible, we say that 
$P$ is \textbf{elliptic}. 
\label{fco.15}\end{definition}

As opposed to the case of usual pseudodifferential operators, ellipticity is 
not a sufficient condition to ensure that an operator is Fredholm.  Some
control near the boundary is also needed.  Let $\fcn Y\cong TY\times\bbR$ 
be the null bundle on $Y$ of the restriction $\fct_{\pa X}X\to TX$.

\begin{proposition}
For any $m\in \bbR$, there is a short exact sequence
\[
 x\fc^{m}(X;E,F)\longrightarrow \fc^{m}(X;E,F) \overset{\no}\longrightarrow
\fs^{m}(\pa X;E,F)
\]
where $\no$ is the \textbf{normal operator} and $\fs^{m}(\pa X;E,F)$ is the
space of $\fcn Y$-suspended fibre pseudodifferential operators on $\pa X$ of 
order
$m$ (see definition $2$ in \cite{mazzeo-melrose4}).  The normal operator
is multiplicative, that is $\no(A)\circ\no(B)=\no(A\circ B)$ whenever $A$ and
$B$ compose. In terms of the
Schwartz kernels, the normal operator $\no$ is just the \textbf{restriction
of the Schwartz kernel to the front face} $\ff_{\Phi}$.
\label{fco.16}\end{proposition}

\begin{definition}
An operator $P\in \fc^{m}(X;E,F)$ is said to be \textbf{fully elliptic} if it 
is elliptic and if $\no(P)$ is invertible.
\label{fco.17}\end{definition}

\begin{definition}
Let $\Ld(X;E)$ be the space of square integrable sections of $E$ with respect
to some metric on $E$ and some smooth positive density on $X$.  For $l\in\bbR$ 
and $m\ge 0$, we define the associated weighted Sobolev spaces by
\begin{equation}
x^{l}\Hc^{m}(X;E)= \{ u\in x^{l}\Ld(X;E)\; | \; Pu\in x^{l}\Ld(X;E)\; \forall
                        P\in \fc^{m}(X;E)\},
\end{equation}
and
\begin{multline}
    x^{l}\Hc^{-m}(X;E)= \\ 
       \left\{ u\in\mathcal{C}^{-\infty}(X;E)\; | \; 
         u=\sum_{i=1}^{N} P_{i}u_{i}, 
        u_{i}\in x^{l}\Ld(X;E),\; P_{i}\in \fc^{m}(X;E) \right\}. 
\end{multline} 
\label{fco.31}\end{definition}
\begin{proposition}
An operator $P\in \fc^{m}(X;E,F)$ is Fredholm as a map
\[
        P:x^{l}\Hc^{m'}(X;E)\to x^{l}\Hc^{m'-m}(X;F) 
\]
for any $l,m\in\bbR$ if and only if it is fully elliptic.
\label{fco.18}\end{proposition}
\begin{proposition}
\textbf{(elliptic regularity)} The null space of a fully elliptic operator
$P\in\fc^{m}(X;E,F)$ is contained in $\CId(X;E)$, the algebra of smooth 
sections of $E$ vanishing in Taylor series at the boundary $\pa X$.
\label{fco.32}\end{proposition}

It would be unfair to end this section without discussing in more details the 
space of suspended operators $\fs^{*}(\pa X;E,F)$.  Since this paper is mostly 
concerned with smoothing operators, we will restrict our attention to 
$\fs^{-\infty}(\pa X;E,F)$.  By taking the \textbf{Fourier transform} in the 
fibres of $\fcn Y$,
it is possible to describe
 $P\in \fs^{-\infty}(\pa X;E,F)$ as a smooth family of 
smoothing operators parametrized by $\fcn^{*}Y$, the dual of $\fcn Y$.  
More precisely, consider the bundle
\begin{equation}
\xymatrix{\Psi^{-\infty}(Z;E,F)\ar@{-}\ar[r]& \mathcal{P}^{-\infty}\ar[d]\\&Y}
\label{fco.19}\end{equation}  
whose fibre at $y\in Y$ is given by $\Psi^{-\infty}(\Phi^{-1}(y);
\left. E \right|_{\Phi^{-1}(y)},\left. F \right|_{\Phi^{-1}(y)})$.
\begin{remark}
When $\dim Z=0$ and $\Phi:\pa X\to Y$ is a finite covering, the bundle 
$\mathcal{P}^{-\infty}$ is the vector bundle $\hom(\mathcal{E},
\mathcal{F})$, where $\mathcal{E}\to Y$ and $\mathcal{F}\to Y$ are the complex
vector bundles over $Y$ with fibres at $y\in Y$ given by
\[
           \mathcal{E}_{y}= \bigoplus_{z\in\Phi^{-1}(y)} E_{z},\quad
         \mathcal{F}_{y}= \bigoplus_{z\in\Phi^{-1}(y)}F_{z}.
\]
\label{fco.30}\end{remark}

We can pull back $\mathcal{P}^{-\infty}$ on $\fcn^{*} Y$ to get a bundle 
$\pi^{*}\mathcal{P}^{-\infty}\to \fcn^{*}Y$, where $\pi: \fcn^{*}Y\to Y$ 
is the 
projection associated to the vector bundle $\fcn^{*}Y$.

\begin{proposition}
There is a one-to-one correspondence between $\fs^{-\infty}(\pa X; E,F)$ and
smooth sections of $\pi^{*}\mathcal{P}^{-\infty}\to \fcn^{*}Y$ which are 
\textbf{rapidly decreasing} with all derivatives as one approaches infinity 
in $\fcn^{*}Y\cong
T^{*}Y\times \bbR$. 
\label{fco.20}\end{proposition}
\begin{remark}
For $P\in \fs^{m}(\pa X; E, F)$, one has a similar correspondence, but one 
must replace $\pi^{*}\mathcal{P}^{-\infty}$ by the appropriate bundle 
$\pi^{*}\mathcal{P}^{m}$ of pseudodifferential operators of order $m$, while 
the sections at infinity should grow no faster than a homogeneous function
of order $m$.  The derivatives of the sections must also satisfy some
growing conditions at infinity.   
\label{fco.24}\end{remark}

This is the point of view we will adopt for the rest of this paper, that is,
we will consider the Fourier transform of $\no(P)$, which is called the 
\textbf{indicial family} of $P$, instead of $\no(P)$.  
One advantage of doing so is that it is relatively
easy to understand the indicial family in terms of local coordinates as
in \eqref{fco.11}.  If $P\in \fc^{-\infty}(X;E,F)$ acts locally as in
proposition~\ref{fco.22} and if $\tau, \eta$ are the dual variables to $x, y$,
then the indicial family is given by
\begin{equation}
\widehat{\no(P)}(y,\eta,\tau)=\int e^{i\eta\cdot Y} e^{i\tau S} 
                           P(0,y,z,S,Y,z-z')dS dY .
\label{fco.25}\end{equation}

\section{The Trace-Defect Formula}\label{tdf.0}

For usual pseudodifferential operators on a closed manifold, it is a well-known
fact that an operator will be of trace class provided its order is 
sufficiently low.  For fibred cusp pseudodifferential operators, this is 
no longer true.  Even smoothing operators are not necessarily of trace
class.  This is because there is a potential singularity at the boundary.

Let us concentrate on the case of smoothing operators.  So we consider 
the space $\fc^{-\infty}(X;E)$, where $E$ is some complex vector 
bundle on $X$.  To avoid possible ambiguities, we assume that in a collar 
neighborhood $\pa X\times [0,1)_{x}\subset X$ of $\pa X$ parametrized by the 
defining function $x$, the vector bundle $E$ is of the form 
$E=\pi^{*}E_{\pa X}$, where $E_{\pa X}$ is a vector bundle on $\pa X$ and 
$\pi: \pa X\times [0,1)_{x}\to \pa X$ is the projection on the first factor.
In other words, we choose an \textbf{explicit trivialization} of $E$ near the
boundary which agrees with the choice of the defining function $x$.  This way,
we can make sense of Taylor series at the boundary in powers of $x$ (see 
\eqref{tdf.3} below).  

If $A\in \fc^{-\infty}(X;E)$ is of trace class,
then its trace is given by integrating its Schwartz kernel along the 
diagonal
\begin{equation}
\Tr(A)=\int_{\Delta_{\Phi}} \tr_{E}K_{A} =\int_{\Delta_{\Phi}} \tr_{E}(K_{A}')
x^{-l-2}\nu, \hspace{1cm}K_{A}'\in\CI(X^{2}_{\Phi};\Hom(E,E)),
\label{tdf.1}\end{equation} 
where $l=\dim Y$, $(x')^{-l-2}\nu$ is a smooth section of the right
fibred cusp density bundle $\fcde_{R}'$, and
$\Hom(E,E)$ is the pull back under the blow-down map 
$\beta_{\Phi}\circ\beta_{b}$ of the homomorphism bundle on $X^{2}$ with 
fibre at $(p,p')\in X^{2}$ given by $\hom(E_{p'},E_{p})$.  When restricted
to the diagonal $\Delta\subset X^{2}$, $\Hom(E,E)$ is the same as $\hom(E,E)$.

The integral \eqref{tdf.1} does not converge for a general 
$A\in\fc^{-\infty}(X;E)$ since $\nu$ is a 
smooth positive
density when restricted to the diagonal $\Delta_{\Phi}$.  To extend
the trace to all smoothing operators, one has to subtract the source of
divergence.  To this end, consider the function
\begin{equation}
        z\mapsto \Tr(x^{z} A),
\label{tdf.2}\end{equation}
which is a priori only defined for $\Re \{z\}>l+1$.  
\begin{lemma}
The function $\Tr(x^{z} A)$ is holomorphic for $\Re \{z\}>l+1$ and it admits
a meromorphic extension to the whole complex plane with at most simple poles
at $l+1-\bbN_{0}$.
\label{tedf.3}\end{lemma}
\begin{proof}
The holomorphy is clear.  Recall from proposition~\ref{fco.16} that in
terms of the Schwartz kernels, the normal operator is the restriction to the 
front face $\ff_{\Phi}$.  Let us consider instead the \textbf{full
normal operator} $\no'$ which gives us the full Taylor series of the
Schwartz kernel at the front face $\ff_{\Phi}$, using $x$ as a defining 
function,
\begin{equation}
   \no'(A)=\sum_{k=0}^{\infty}x^{k}A_{k}, \hspace{1cm} A_{k}\in 
     \fs^{-\infty}(\pa X;E).
\label{tdf.3}\end{equation}
This Taylor series does not necessarily converge, so $\no'(A)$ is really
a formal power series in $x$ which contains all the information about the
derivatives of $A$ at the front face $\ff_{\Phi}$.
The Taylor series $\no'(A)$ contains all the terms causing a divergence in the 
definition of $\Tr(x^{z} A)$.  A simple computation shows that the function
$\Tr(x^{z}x^{k}A_{k}\chi(x))$, which is a priori only defined for
$\Re\{z\}>l+1-k$, has a meromorphic extension to the whole complex plane with a
single simple pole at $z=l+1-k$.  Here, $\chi\in\CI_{c}([0,+\infty))$ is some
cut-off function with $\chi\equiv 1$ in a neighborhood of $x=0$.
Thus, we see that $\Tr(x^{z} A)$ admits a meromorhic extension with at most
simple poles at $l+1-\bbN_{0}$.
\end{proof}

In particular, there is in general a pole at $z=0$ coming from the term
$A_{l+1}$ in \eqref{tdf.3}.
\begin{definition}
For $A\in \fc^{-\infty}(X;E)$, the \textbf{residue trace} $\rTr(A)$ of 
$A$ is the residue at $z=0$ of the meromorphic function $z\mapsto 
\Tr(x^{z} A)$.
\label{tdf.4}\end{definition}
Using the local representation of proposition~\ref{fco.22}, the residue trace
 can be expressed as
\begin{equation}
\rTr(A)=\int \tr_{E}(A_{l+1}(y,z,S=0,Y=0,z-z'=0))dy dz.
\label{tdf.5}\end{equation}
Taking the Fourier transform in the fibres of $\no Y\cong TY\times \bbR$ 
as in \eqref{fco.25}, the residue trace can be rewritten as (cf. equations
4.1 and 4.2 in \cite{mr96h:58169})
\begin{equation}
\rTr(A)= \frac{1}{(2\pi)^{l+1}}\int_{\no^{*}Y}  \Tr_{\Phi^{-1}(y)}
(\widehat{A}_{l+1}(y,\tau,\eta)) dy d\eta d\tau
\label{tdf.6}\end{equation}
where $\Tr_{\Phi^{-1}(y)}$ is the trace on $\Psi^{-\infty}(\Phi^{-1}(y);
\left. E\right|_{\Phi^{-1}(y)})$.  Notice that this integral does not
depend on the choice of coordinates on $Y$, since $dy d\eta$ is the volume
form coming from the canonical symplectic form on $T^{*}Y$.  Of course,
$d\tau$ depends on the choice of the defining function $x$, but we assume
it has been fixed.

One can then extend the trace to all smoothing operators in 
$\fc^{-\infty}(X;E)$ by subtracting the pole at $z=0$.

\begin{definition}
The \textbf{regularized trace} $\bTr(A)$ of $A\in \fc^{-\infty}(X;E)$ is given
by
\[
    \bTr(A)= \lim_{z\to 0} \left( \Tr(x^{z}A)- \frac{\rTr(A)}{z}\right) .
\]
\label{tdf.7}\end{definition}

\begin{remark}
For $A\in x^{l+2}\fc^{-\infty}(X;E)$, $A$ is of trace class and
$\bTr(A)=\Tr(A)$.
\label{tdf.12}\end{remark}

The regularized trace is defined for all $A\in \fc^{-\infty}(X;E)$ and it
is an extension of the usual trace, but it is not properly speaking a trace,
since it does not vanish on commutators.  The trace-defect formula measures
how far the regularized trace is from being a trace.

\begin{proposition}
\textbf{(Trace-defect formula)} For $A,B\in \fc^{-\infty}(X;E)$, the 
trace-defect is given by
\[
     \bTr([A,B])= \rTr( (D_{\log x}A)B)
\]
where $D_{\log x}A=\left. \frac{\pa}{\pa z}A_{z}\right|_{z=0}$ and 
$A_{z}=x^{-z}[x^{z},A]$.
\label{tdf.8}\end{proposition}
\begin{proof}
For $\Re \{z\}\gg 0$, we have
\begin{equation}
\begin{aligned}
\Tr(x^{z}[A,B]) &=\Tr(x^{z}AB- x^{z}BA)=\Tr(x^{z}AB) -\Tr((x^{z}B)A) \\
        &= \Tr(x^{z}AB) -\Tr(Ax^{z}B) -\Tr([x^{z}B,A]).
\end{aligned}
\label{tedf.9}\end{equation}
But $\Tr([x^{z}B,A])= 0$ for $\Re \{z\}\gg 0$, so
\begin{equation}
\begin{aligned}
\Tr(x^{z}[A,B]) &=\Tr(x^{z}AB) -\Tr(Ax^{z}B) \\
        &=\Tr([x^{z},A]B)=\Tr(x^{z}A_{z}B). 
\end{aligned}
\label{tdf.10}\end{equation}
Here, $A_{z}=x^{-z}[x^{z},A]$ is entire in $z$ and vanishes at $z=0$, so
\begin{equation}
\Tr(x^{z}[A,B])= z\Tr(x^{z}(D_{\log x}A) B) + \mathcal{O}(z), \hspace{0.5cm}
     \mbox{where} \hspace{0.3cm}D_{\log x}A=
\left. \frac{\pa}{\pa z}A_{z} \right|_{z=0}.
\label{tdf.11}\end{equation}
Thus, we see that $\Tr(x^{z}[A,B])$ has no pole at $z=0$ and its value at 
$z=0$ is given by
\[
          \bTr([A,B])= \rTr( (D_{\log x}A) B ).
\]
\end{proof}

In the case of cusp operators ($Y=\{\pt\}$), it is possible to use 
\eqref{tdf.6} to get an explicit formula depending only on the 
indicial families of $A$ and $B$ (cf. \cite{fipomb}).  In the more general
case, higher order terms of the Taylor series \eqref{tdf.3} are involved.
Combined with the geometry of $Y$, this makes an explicit computation
harder or maybe even impossible to get.

\section{The Index in terms of the Trace-Defect Formula}\label{itdf.0}

Let $\Frc(X;E)$ denote the space of Fredholm operators of the form 
$\Id + A$ with $A\in \fc^{-\infty}(X;E)$.  From proposition~\ref{fco.18},
the space $\Frc(X;E)$ is given by
\begin{equation}
\Frc(X;E)=\{\Id +A \;\; |\; \; A\in \fc^{-\infty}(X;E),\;\; 
\no(\Id+A) \hspace{0.2cm}\mbox{is invertible}\}.
\label{itdf.1}\end{equation}

The goal of the next few sections will be to obtain a topological description
of the index of operators in $\Frc(X;E)$, that is, a description in terms of
$K$-theory.  As a first step, let us use the trace-defect formula to get some
insight about the index.

\begin{notation}
From now on we will follow the notation of \eqref{tdf.3} for the indicial
family.  So 
$\widehat{A}_{0}=\widehat{\no(A)}$ will denote  
the indicial family.  In proposition~\ref{fco.20}, it is $\widehat{A}_{0}$
which can be seen as a section of the bundle $\pi^{*}\mathcal{P}^{-\infty}$. 
\label{itdf.2}\end{notation}

Given $(\Id+A)\in \Frc(X;E)$, let $(\Id+B)\in\Frc(X;E)$ be a parametrix, that
is
\begin{equation}
  (\Id+A)(\Id+B)=\Id +Q_{1},\hspace{0.5cm} (\Id+B)(\Id+A)=\Id +Q_{2},
\label{itdf.3}\end{equation}
with $Q_{1},Q_{2}\in x^{\infty}\fc^{-\infty}(X;E)$, where 
\begin{equation}
x^{\infty}\fc^{-\infty}(X;E)=\{ Q\in \fc^{-\infty}(X;E)\;\; | \;\;
                                   \no'(Q)=\sum_{k=0}^{\infty}x^{k}Q_{k}=0 \}. 
\label{itdf.4}\end{equation}
In particular, $Q_{1}$ and $Q_{2}$ are compact operators of trace class,
so $[\Id+A,\Id+B]=Q_{1}-Q_{2}$ is also of trace class.  Using Calder\'{o}n's 
formula for the index, we have that
\begin{equation}
\begin{aligned}
\ind (\Id+A)&= \Tr([\Id+A,\Id +B]) \\
            &= \bTr([\Id+A,\Id +B]), \hspace{0.3cm}
                                    \mbox{by remark~\ref{tdf.12}}\\
            &=\bTr([A,B]) \\
            &=\rTr((D_{\log x}A) B), \hspace{0.3cm}\mbox{by 
                                              proposition~\ref{tdf.8}.}
\end{aligned}
\label{itdf.5}\end{equation}

Using \eqref{tdf.6}, we get the following formula for the index.

\begin{proposition}
The index of $(\Id +A)\in \Frc(X;E)$ is given by
\[
           \ind(\Id+A) = \frac{1}{(2\pi)^{l+1}}\int_{\fcn^{*}Y}  
\Tr_{\Phi^{-1}(y)}
(\widehat{(D_{\log x}A B)}_{l+1}(y,\tau,\eta)) dy d\eta d\tau
\]
where $(\Id+B)\in \Frc(X;E)$ is a parametrix for $A$.
\label{itdf.6}\end{proposition}

In the case of cusp operators, one can rewrite this formula in terms
of the indicial
family of $A$ only (see \cite{fipomb}).  In fact, one can show without 
using the trace-defect formula that the index only depends on the indicial
family.

\begin{proposition}
The index of $(\Id+A)\in \Frc(X;E)$ only depends on the homotopy class
of the indicial family $\Id +\widehat{A}_{0}$
seen as an \textbf{invertible} section of the bundle 
$\pi^{*}\mathcal{P}^{-\infty}$ of
proposition~\ref{fco.20} which \textbf{converges rapidly} to $\Id$ with all 
derivatives as one approaches infinity.
\label{itdf.7}\end{proposition}
\begin{proof}
Suppose that $t\mapsto \Id + \widehat{A}_{0}(t), \, t\in[0,1]$ is a
(smooth) homotopy of invertible indicial families with $\widehat{A}_{0}(t)\in 
\fs^{-\infty}(\pa X;E)$ for all $t\in[0,1]$.  Suppose also that 
$\Id+ \widehat{A}_{0}(0)$ and $\Id+ \widehat{A}_{0}(1)$ are the indicial 
families of $(\Id+A(0)), (\Id+A(1))\in \Frc(X;E)$.  Thinking in terms of
Schwartz kernels and using Seeley extension for manifolds with corners (see
the proof of propostion 1.4.1 in \cite{MelroseMWC}), we can construct a 
homotopy $t\mapsto A(t), t\in[0,1]$, joining
$A(0)$ and $A(1)$, such that $A(t)\in\fc^{-\infty}(X;E)$ and 
\[
           \widehat{\no(A(t))}=\widehat{A}_{0}(t) \hspace{0.2cm}\forall 
     t\in[0,1].
\]
Since $t\mapsto \Id+\widehat{A}_{0}(t)$ is a homotopy of \textbf{invertible}
indicial families, we know by proposition~\ref{fco.18} that 
$t\mapsto(\Id +A(t))$ is a homotopy of Fredholm operators.  By the homotopy
invariance of the index, we conclude that 
\[
                 \ind(\Id +A(0))=\ind(\Id +A(1)),
\] 
which shows that the index of an operator $(\Id+A)\in\Frc(X;E)$ only depends on
the homotopy class of its indicial family.
\end{proof}

One therefore expects that it is possible, at least in theory, to rewrite the
formula of proposition~\ref{itdf.6} so that it only involves the indicial
family.  A first step in this direction is to notice that, by 
proposition~\ref{itdf.7}, we can assume that the operator $\Id+A\in\Frc(X;E)$
is such that 
\[
           \widehat{\no'(A)}= \widehat{A}_{0}, \; \mbox{with}\; A_{k}=0\;
                \forall k\in \bbN.
\]
However, once this assumption is made on $A$, it is not possible to assume the
same for the parametrix $B$.  To be more precise, the equation
\eqref{itdf.3} completely determined $\widehat{\no'(B)}$ in terms of
$\widehat{\no'(A)}$.  This will become clear through the concept of 
$*$-product.
\begin{definition}
Let $\fs^{-\infty}(\pa X;E)[[x]]$ denote the space of formal power series in
$x$ with coefficients taking values in $\fs^{-\infty}(\pa X;E)$.
\label{itdf.8}\end{definition}
\begin{lemma}
The full normal operator $\no': \fc^{-\infty}(X;E)\to 
\fs^{-\infty}(\pa X; E)[[x]]$ is surjective.
\label{itdf.9}\end{lemma}
\begin{proof}
This is an easy consequence of Seeley extension for manifolds with corners.
\end{proof}
\begin{definition}
For $a,b\in \fs^{-\infty}(\pa X; E)[[x]]$, the $*$-product 
$a*b\in \fs^{-\infty}(\pa X; E)[[x]]$ is defined to be
\[
   a*b= \widehat{\no'(A\circ B)}
\] 
where $A, B\in\fc^{-\infty}(X;E)$ are chosen such that $\widehat{\no'(A)}=a$, 
$\widehat{\no'(B)}=b$.
It does not depend on the choice of $A$ and $B$, since 
$x^{\infty}\fc^{-\infty}(X;E)$ is an ideal of $\fc^{-\infty}(X;E)$.
\label{itdf.10}\end{definition}
\begin{proposition}
There exist differential operators $P_{k,j,l}$ and $Q_{k,j,l}$ 
acting on sections
of the bundle $\mathcal{P}^{-\infty}$ defined in \eqref{fco.19} so that,
for any $a,b\in \fs^{-\infty}(\pa X; E)[[x]]$, the $*$-product of $a$ and 
$b$ is given by
\[
  a*b=\left( \sum_{j=0}^{\infty}x^{j}a_{j}\right)*
              \left( \sum_{j=0}^{\infty}x^{j}b_{j}\right)=
   \sum_{k=1}^{\infty}\;\sum_{j+l\le k}  x^{k} (P_{k,j,l}a_{j})\circ
      (Q_{k,j,l}b_{l}).
\] 
Moreover, if $\Phi:\pa X\to Y$ is a trivial fibration, then, once
an explicit trivialization has been chosen, the differential operators
$P_{k,j,l}$ and $Q_{k,j,l}$ are independent of the geometry of the typical 
fibre $Z$.
\label{itdf.11}\end{proposition}
For the proof of proposition~\ref{itdf.11}, we refer to proposition $3.11$
in \cite{Lauter-Moroianu2}.  In terms of the $*$-product, 
equation \eqref{itdf.3} translates into
\begin{equation}
    \widehat{\no'(A)} + \widehat{\no'(B)} =-\widehat{\no'(A)}*
       \widehat{\no'(B)}= -\widehat{\no'(B)}*\widehat{\no'(A)}.
\label{itdf.20}\end{equation}
Looking at this equation order by order in $x$ completely determines
$\widehat{\no'(B)}$ in terms of $\widehat{\no'(A)}$.  In general, for 
$k\in\bbN$, on can expect 
that $\widehat{B}_{k}\ne 0$, so $\widehat{\no'(B)}$ will  involve a term of 
order $x^{k}$.
This is also the case for $D_{\log x}A$.

\begin{lemma}
If $\widehat{\no'(A)}=\widehat{A}_{0}$, then 
\[
\widehat{\no'(D_{\log x}A)}=
\sum_{k=1}^{\infty} \frac{x^{k}D_{\tau}^{k}\widehat{A}_{0}}{k}.
\]
\label{itdf.12}\end{lemma}
\begin{proof}
The Schwartz kernel of $A_{z}=x^{-z}[x^{z},A]$ is 
$A(1-\left(\frac{x'}{x}\right)^{z})$, so the the Schwartz kernel of 
$D_{\log x} A$ is given by
\[
      \frac{\pa}{\pa z} \left.\left( A(1-\left(\frac{x'}{x}\right)^{z}) \right)
\right|_{z=0}= -A\log\left( \frac{x'}{x}\right).
\]
But $\frac{x'}{x}=1-Sx$, where $S=\frac{x-x'}{x^{2}}$, so
\[
   D_{\log x}A = -A \log(1-Sx) = A \sum_{k=1}^{\infty} \frac{(Sx)^{k}}{k}.
\]
From \eqref{fco.25}, $\widehat{\no(S^{k}A)}=D_{\tau}^{k}\widehat{A}_{0}$,
where $D_{\tau}=-i\frac{\pa}{\pa \tau}$, so we conclude that
\[
 \widehat{\no'(D_{\log x}A)}= \sum_{k=1}^{\infty} \frac{x^{k}D_{\tau}^{k}
\widehat{A}_{0}}{k}.
\]
\end{proof}

Thus, knowing the differential operators $P_{k,j,l}$ and $Q_{k,j,l}$, we
see that it is possible to rewrite the formula of proposition~\ref{itdf.6}
into an expression involving only the indicial family of $A$.  The problem
of course is to actually find out what are those differential operators.
This is relatively easy in the cusp case, but the case of an arbirary 
manifold $Y$ seems hopeless.
We will not try to develop further the formula of proposition~\ref{itdf.6}.  As
it is written now and from proposition~\ref{itdf.11}, one can extract some 
important information about the index.
For instance, notice that the formula does not involve the geometry of the
interior of $X$, so if $X'$ is another manifold with the same boundary as 
$X$ and if $(\Id+A')\in \Frc(X';E)$ has the same indicial family as the one
for $\Id+A$, then $\Id+A$ and $\Id + A'$ have the same index.  In 
section~\ref{im.0},
the formula of proposition~\ref{itdf.6} will play a crucial r\^{o}le in 
reducing the computation of the index to the case of a scattering operator.  
For the moment however, we will try to have a better understanding of the
topological information contained in the indicial family.
When the fibration $\Phi$ is trivial, it is almost straightforward to define
a $K$-class out of the homotopy class of the indicial family.  But to include 
the more general case of a non-trivial fibration, we first need to discuss the 
concept of spectral section.

\section{Spectral Sections}
\label{ss.0}

Spectral sections were originally introduced in 
\cite{melrose-piazza} to describe the boundary conditions for families 
of Dirac operators on manifolds with boundary.  We intend to use spectral 
sections for a very different purpose, and instead of dealing with families
of Dirac operators, we will consider families of Laplacians.

Let $\phi:M\to B$ be a smooth (locally trivial) fibration of compact manifolds,
where the fibres are closed manifolds and the base $B$ is a compact
manifold with possibly a non-empty boundary $\pa B$.  
Let $T(M/B)\subset TM$ denote the null space of the differential
\[\phi_{*}:TM\to TB\] 
of $\phi$.  Clearly, $T(M/B)$ is a subbundle of $TM$ which
on each fibre $M_{b}=\phi^{-1}(b)$  restricts to be canonically 
isomorphic to $TM_{b}$.  Let $g_{M/B}$ be a family of metrics on $M$, that is,
$g_{M/B}$ is a metric on $T(M/B)$ which gives rise to a Riemannian metric on
each fibre.  Let $E$ be a complex vector bundle with an Euclidean metric 
$g^{E}$ and a connection $\nabla^{E}$.

This allows us to define a smooth family of Laplacians $\Delta_{M/B}$, which
on a fibre $M_{b}=\phi^{-1}(b)$ acts as
\begin{equation}
\Delta_{b}s= -\Tr( \nabla^{T^{*}M_{b}\otimes E}\nabla^{E}s), \hspace{0.2cm}
   s\in \CI(M_{b};E_{b}), \hspace{0.2cm}E_{b}=\left. E\right|_{M_{b}},
\label{ss.1}\end{equation} 
where $\nabla^{T^{*}M_{b}\otimes E}$ is induced from $\nabla^{E}$ and the
Levi-Civita connection of $g_{b}=\left. g_{M/B}\right|_{M_{b}}$, and
$\Tr(S)\in \CI(M_{b};E_{b})$ is the contraction of an element 
\[S\in\CI(M_{b}; T^{*}M_{b}\otimes T^{*}M_{b}\otimes E_{b})\]
with the
metric $g_{b}\in \CI(M_{b}, TM_{b}\otimes TM_{b})$. 

Let $\Ld(M/B;E)\to B$ be the Hilbert bundle over $B$ with fibre at $b\in B$ 
given by $\Ld(M_{b};E_{b})$.  The scalar product on $\Ld(M_{b};E_{b})$ is 
defined in the usual way by
\begin{equation}
  \langle s_{1}, s_{2}\rangle= \int g^{E}(s_{1},\overline{s}_{2})\;dg_{b}, 
\hspace{0.4cm} s_{1},s_{2}\in\CI(M_{b};E_{b}),
\label{ss.2}\end{equation} 
where $dg_{b}$ is the volume form associated to the metric $g_{b}$.

As is well-known, the Laplacian $\Delta_{b}$ is self-adjoint and has a
non-negative discrete spectrum.  Moreover, eigensections with different
eigenvalues are orthogonal in $\Ld(M_{b};E_{b})$.

\begin{definition}
A \textbf{spectral section} $(\Pi_{M/B},R_{1},R_{2})$ for the family of 
Laplacians $\Delta_{M/B}$ is a 
smooth family $\Pi_{M/B}:\Ld(M/B;E)\to \Ld(M/B;E)$ of projections
\[ 
\Pi_{b}:\Ld(M_{b};E_{b})\to \Ld(M_{b};E_{b}), \hspace{0.3cm}b\in B,
\]
with range a \textbf{trivial} vector bundle over $B$, together with real 
numbers $R_{1}$ and $R_{2}$, $R_{1}<R_{2}$, such that for all $b\in B$
\[
    \Delta_{b}s=\lambda s \hspace{0.3cm}\Longrightarrow \hspace{0.2cm}
\left\{ \begin{array}{r@{\quad \mbox{if}\quad}l}
         \Pi_{b}s=s & \lambda< R_{1}, \\
         \Pi_{b}s=0 & \lambda>R_{2}.
        \end{array}\right.
\]
\label{ss.3}\end{definition}
\begin{remark}
The triviality of the range of the spectral section  is not in the original
definition.  We added this property because we will only consider spectral
sections having as range  a trivial vector bundle over $B$.  Notice also that
from the discreteness of the spetrum, this vector bundle must be of finite 
rank.
\label{ss.6}\end{remark}

Obviously, one gets a spectral section by taking $\Pi_{M/B}=0$ and 
$R_{1}<R_{2}<0$, but this example is rather trivial.  We would like to 
know about the existence of non-trivial spectral sections.  The following 
result, which was communicated to the author in a private discussion,
 is due to R.B. Melrose.

\begin{proposition}
Let $\Delta_{M/B}$ be a family of Laplacians as in \eqref{ss.1}.  Then given
$R>0$, there exists $R'>R$ and a smooth family of projections
$\Pi_{M/B}$ such that $(\Pi_{M/B},R,R')$ is a spectral section for
$\Delta_{M/B}$.
\label{ss.4}\end{proposition}
\begin{proof}
By compactness, the discreteness of the spectrum and its continuity as a 
set-valued function, we can find a covering of $B$ by a finite number 
of open sets $\Omega_{i}$, and $R_{i}\in(R,\infty)$ such that $R_{i}$ is
not in the spectrum of $\Delta_{b}$ for all $b\in \Omega_{i}$.  The span of the
eigensections with eigenvalues less than $R_{i}$ is then a smooth vector
bundle $E_{i}$ over $\Omega_{i}$.  Moreover, there are smooth bundle 
inclusions on all non-trivial intersections 
$\Omega_{ij}=\Omega_{i}\cap\Omega_{j}$, $I_{ij}:E_{i}\to E_{j}$, provided 
$R_{i}\le R_{j}$.  These inclusions satisfy obvious compatibility conditions
on triple intersections.

By Kuiper's theorem (see \cite{Kuiper}), the Hilbert bundle $\Ld(M/B;E)$ is
trivial in the uniform topology.  It follows that there are trivial 
finite dimensional subbundles of arbitrary large rank. In other words, there 
is a sequence of families of
finite rank projections $\pi^{(k)}_{b}$ with ranges  
 trivial vector bundles over $B$  such that $\pi_{b}^{(k)}\to
\Id$ in the strong topology (on operators).  
Setting $\pi_{b}=\pi_{b}^{(k)}$ for sufficiently 
large $k$, it follows that the restriction of $\pi_{b}$ to the $E_{i}$ are 
injective
\[
           \pi_{b}:E_{i,b}\hookrightarrow F_{b},
\] 
where the trivial vector bundle $F$ is the range of $\pi$.
By taking a larger $k$ if needed, we can also assume that for all $E_{i}$, the
norm of $\left. (\Id-\pi_{b})\right|_{E_{i}}$ is less than $\frac{1}{2}$.  
These embeddings of the bundles $E_{i}$ as subbundles of 
$\left. F\right|_{\Omega_{i}}$ are consistent with the inclusions $I_{ij}$.  
We may also consider the generalized inverse of $\pi_{b}$ on $E_{i}$,
 $m_{i}:F\to E_{i}$ over $\Omega_{i}$, defined as the composite of the 
orthogonal projection on $\pi(E_{i})$ and the inverse of $\pi$ as a map
from $E_{i}$ to $\pi(E_{i})$.  

Let $\mu_{i}$ be a partition of unity subordinate to the $\Omega_{i}$, and 
consider the family of linear maps 
\begin{equation}
  f_{b}:F_{b}\to \Ld(M_{b};E_{b}),\hspace{0.5cm} 
  f_{b}=(\Id-\pi_{b})\sum_{j}\mu_{j}(b)m_{j}(b).
\label{ss.5}\end{equation} 
This is a well-defined smooth family.  Consider the open sets 
\[
     \mathcal{U}_{i}=\Omega_{i}\setminus \bigcup_{\{j; R_{j}<R_{i}\} } 
     \supp(\mu_{j}).
\]
These cover $B$.  Over $\mathcal{U}_{i}$, the sum in \eqref{ss.5} is
limited to those $j$ with $\Omega_{ij}\ne\emptyset $ and 
$m_{j}(\pi_{b}e)=e$ for $e\in E_{i,b}$.  Therefore, $f_{b}(\pi_{b}e)=
(\Id-\pi_{b})e$ for $b\in \mathcal{U}_{i}$ and $e\in E_{i,b}$.  Since
$\| f_{b}\| < \frac{1}{2}$, it follows that 
\[
          (\Id+f_{b}): F_{b}\to \Ld(M_{b};E_{b})
\]
embeds $F$ as a subbundle of $\Ld(M/B;E)$ which, over $\mathcal{U}_{i}$,
contains $E_{i}$ as a subbundle.

Finally, choose $\psi\in \CI(\bbR)$ with $0\le \psi(t)\le 1$ and such that
\[
          \psi(t)= \left\{ \begin{array}{r@{\quad \mbox{if}\quad}l}
                                 1 & t<0 \\
                                 0 & t>1.
                            \end{array} \right.
\]
For $T\in\bbR $, consider the family of linear maps
$Q_{T}=\psi(\Delta_{M/B}-T)$.  For large $T$,
$Q_{T}\circ(\Id+f):F\to G$ embeds $F$ as a spectrally finite subbundle $G$ 
containing the $E_{i}$ over $\mathcal{U}_{i}$.  Then, if $\Pi_{M/B}$ is
the orthogonal projection on $G$ and $R'=T+1$, we see that $(\Pi_{M/B},R,R')$
is the desired spectral section.  Since $F$ is a trivial vector
bundle, it is clear that $G$ is a trivial vector bundle as well.
\end{proof}

\section{Construction of the $K$-Class}\label{ckc.0}

We are now ready to construct a $K$-class out of the indicial family of a 
Fredholm operator in $\Frc(X;E)$.  We need to distinguish two cases, the
case where $\dim Z>0$ and the case where $\dim Z=0$.

\subsection{The case where $\dim Z>0$}
\label{ckc.1}

Let
\begin{equation}
\xymatrix{\G(Z;E)\ar@{-}\ar[r]&\G\ar[d]^{\phi}\\&Y}
\label{ckc.2}\end{equation}
be the smooth bundle with fibre at $y\in Y$ given by $\G(Z_{y};E_{y})$,
where
\begin{equation}
  \G(Z_{y};E_{y})=\{ \Id +Q \; \; | \;\; Q\in\Psi^{-\infty}(Z_{y};E_{y}),
  \; \Id +Q \; \mbox{is invertible} \}
\label{ckc.3}\end{equation}
is the group of invertible smoothing perturbations of the identity discussed
in the introduction and 
$Z_{y}=\Phi^{-1}(y)$, $E_{y}=\left. E\right|_{\Phi^{-1}(y)}$.
The bundle $\G$ is a subbundle of the bundle $\mathcal{P}^{-\infty}$ in 
\eqref{fco.19}.
Let $\pi:V\to Y$ be a real vector bundle over $\bbR$.  In many of the 
situations interesting us, the vector bundle $V$ will be $\fcn^{*}Y$.
Let $\pi^{*}\G$ denote the pull-back of $\G$ on $V$. 
\begin{definition}
A smooth section of $\pi^{*}\G$, or more generally of any bundle which has 
a section $\Id$, is said to be \textbf{asymptotic to the identity}
if it converges rapidly with all derivatives to $\Id$ as one approaches 
infinity in $V$.  Let $\Gamma_{\Id}(V;\pi^{*}\G)$ denote the space of 
smooth sections of $\pi^{*}\G$ asymptotic to the identity. 
\label{ckc.1.1}\end{definition}
The following
corollary is an immediate consequence of proposition~\ref{fco.18} and 
proposition~\ref{fco.20}.

\begin{corollary}
For $V=\fcn^{*}Y$, the indicial family $\Id +\widehat{A}_{0}$ of a Fredholm 
operator
in $\Frc(X;E)$ corresponds to a section of the bundle $\pi^{*}\G\to \fcn^{*}Y$ 
which is asymptotic to the identity.
\label{ckc.4}\end{corollary}

Our primary goal in this section is to construct a $K$-class out of the 
\textbf{homotopy class} of the
indicial family of a Fredholm operator $(\Id +A)\in \Frc(X;E)$.
But considering an arbitrary real vector bundle $V$ over $Y$ will allow
us to get a more general result which will turn out to be useful in 
section~\ref{hg.0}.  So we will consider the more general situation of 
a smooth section $\Id+S$ of $\pi^{*}\G\to V$  asymptotic to the identity, 
bearing in mind that in the case where $V=\fcn^{*}Y$, the section
$\Id+S$ can be
seen as the indicial family of a Fredholm operator in $\Frc(X;E)$.

We want to describe the homotopy classes of such sections.   
To this end,
choose a family of metrics $g_{\pa X /Y}$ for the fibration $\Phi:\pa X\to Y$
and a Euclidean metric $g^{E}$ as well as a connection $\nabla^{E}$ for the 
bundle $E$.
Then, as discussed in section~\ref{ss.0}, there is an associated family of 
Laplacians $\Delta_{\pa X/ Y}$.  Given $y\in Y$, recall from the introduction 
that if $\{f_{i}\}_{i\in\bbN}$ is an orthonormal basis of eigensections of 
$\Delta_{y}$ with increasing eigenvalues, then
\begin{equation}
\begin{array}{rrcl}
F_{ij}: & \G(Z_{y};E_{y}) & \rightarrow & \mathcal{G}^{-\infty} \\
        & \Id + Q & \mapsto & \delta_{ij} +\langle f_{i}\; ,\;
Qf_{j}\rangle_{\Ld} 
\end{array}
\label{ckc.5}\end{equation}
is an isomorphism between $\G(Z_{y};E_{y})$ and the group 
$\mathcal{G}^{-\infty}$ of invertible semi-infinite matrices 
$\delta_{ij}+ A_{ij}$ such that 
\begin{equation}
  \sup_{i,j}(i+j)^{k}|A_{ij}|<\infty, \;\; \forall k\in\bbN_{0},
\label{ckc.6}\end{equation}
where $\delta_{ij}$ corresponds to the semi-infinite identity matrix.  If 
$P_{y}(T)$ denotes the projection onto the span of the eigensections
of $\Delta_{y}$ with eigenvalues less than $T$, this means that 
\begin{equation}
  P_{y}(T)(Q_{y})P_{y}(T)\underset{T\to+\infty}{\longrightarrow} Q
\label{ckc.7}\end{equation}
in the $\CI$-topology, so we can approximate $Q$ by spectrally finite 
smoothing operators.

\begin{lemma}
Let $(\Id+S)$ be a smooth section of $\pi^{*}G^{-\infty}\to V$ asymptotic
to the identity, where $V\to Y$ is a real vector bundle over $Y$.  Then there 
exists a 
spectral section $(\Pi, R_{1},R_{2})$ for $\Delta_{\pa X/Y}$ such that 
\[
       \| S_{b}-\Pi_{b}S_{b}\Pi_{b} \| < \frac{1}{2\| (\Id +S_{b})^{-1}\|},
\quad \forall b\in V,          
\] 
where $\|\cdot\|$ is the operator norm and $\Pi_{b}$ denotes $\Pi_{\pi(b)}$. 
\label{ckc.8}\end{lemma}
\begin{proof}
Using \eqref{ckc.7} and the fact $\Id +S$ is asymptotic to the identity, we see
by the compactness of $Y$, the discreteness
of the spectrum and its continuity as a set-valued function, that there
exist a covering of $Y$ by a finite number of open sets $\Omega_{i}$, and 
numbers $T_{i}\in\bbR^{+}$ such that
\begin{equation}
   \|S_{b}-P_{b}(T_{i})S_{b}P_{b}(T_{i})\| <
                         \frac{1}{2\|(\Id+S_{b})^{-1}\|}\;,
\quad \forall\, b\in\pi^{*}\Omega_{i},
\label{ckc.9}\end{equation}
where $P_{b}(T_{i})$ denotes $P_{\pi(b)}(T_{i})$.  
Set $R_{1}=\underset{i}{\max}\{T_{i}\}$.  Then, by proposition~\ref{ss.4},
there exists a smooth family of projections $\Pi$ and 
$R_{2}>R_{1}$ such that $(\Pi,R_{1},R_{2})$ is a spectral section
for $\Delta_{\pa X/Y}$.  By construction, we have that for $b\in
\pi^{*}\Omega_{i}$,
\[
                     \Pi_{b}P_{b}(T_{i})\Pi_{b}=P_{b}(T_{i}),
\]
so the estimate \eqref{ckc.9} also holds for $\Pi_{b}$ for all $b\in 
V$.  Thus $(\Pi,R_{1},R_{2})$ is the desired spectral section.
\end{proof}

\begin{corollary}
Let $(\Pi,R_{1},R_{2})$ be the spectral section of lemma~\ref{ckc.8}.  Then,
$\Id+S$  and $\Id + \Pi S\Pi$ are homotopic as sections of $\pi^{*}\G$ 
asymptotic to the identity.
\label{ckc.10}\end{corollary}
\begin{proof}
Consider the smooth homotopy 
\begin{equation}
  t\mapsto \Id + S +t(\Pi S\Pi-S), \quad t\in [0,1]. 
\label{ckc.11}\end{equation}
For all $b\in V$ and $t\in [0,1]$, notice that
\[
    (\Id+ S_{b})^{-1}(\Id +S_{b}+t(\Pi_{b}S_{b}\Pi_{b}-S_{b}))=
                 \Id +t(\Id+S_{b})^{-1}(\Pi_{b}S_{b}\Pi_{b}-S_{b})
\]
is invertible, since $\|t(\Id+S_{b})^{-1}(\Pi_{b}S_{b}\Pi_{b}-S_{b})\|
<\frac{1}{2}$.
This means that 
\[ \Id +S_{b}+t(\Pi_{b}S_{b}\Pi_{b}-S_{b}) \] 
is invertible as well,
and so the homotopy \eqref{ckc.11} between $(\Id+S)$ and $(\Id+\Pi S\Pi)$
is through sections of $\pi^{*}\G$ asymptotic to the identity.
\end{proof}

Let $F$ denote the range of $\Pi$.  Recall by our definition of a spectral 
section that $F$ is a trivial vector bundle of finite rank over $V$.
Let $\GL(F,\bbC)\to \fcn^{*}Y$ be the smooth bundle over $V$ with 
fibre at $b\in V$ given by $\GL(F_{\pi(b)},\bbC)$, the group of 
complex linear isomorphism of $F_{\pi(b)}$.  Then corollary~\ref{ckc.10} 
effectively
reduces the section $\Id+S$ to a section
\begin{equation}
   \Id+S_{\Pi}=\Pi(\Id+S)\Pi: V \to \GL(F,\bbC).
\label{ckc.12}\end{equation} 
If we also choose an explicit trivialization of $F$ and if $n$ is the rank of 
$F$, then the section
$\Id+S$ can be seen as smooth map
\begin{equation}
   \Id+S_{\Pi}: V \to \GL(n,\bbC)
\label{ckc.13}\end{equation} 
asymptotic to the identity. 
Finally, since $\GL(n,\bbC)$ is a topological subspace of the direct limit
\begin{equation}
        \GL(\infty,\bbC)=\lim_{k\to\infty}\GL(k,\bbC),
\label{ckc.14}\end{equation}
we can think of \eqref{ckc.13} as a map
\begin{equation}
   \Id+S_{\Pi}: V \to \GL(\infty,\bbC)
\label{ckc.15}\end{equation} 
which converges to the identity at infinity.  We are mostly interested
in the homotopy class of this map.  

\begin{definition}
Let $[V \; ;\; \GL(\infty,\bbC)]$ denote the homotopy classes 
of continuous maps from $V$ to $\GL(\infty,\bbC)$ which converges
to the identity as one approaches infinity in $V$.  If $a$ is
such a continuous map, then let [a] denote its homotopy class in 
$[ V \; ;\; \GL(\infty,\bbC)]$.
\label{ckc.16}\end{definition}
\begin{remark}
The passage from the category of smooth maps to the category of 
continuous maps does not change the set of homotopy classes, see 
for instance proposition $17.8$ in \cite{Bott-Tu}.
\label{ckc.17}\end{remark}
\begin{definition}
If $(\Id+S)$ and $(\Pi,R_{1},R_{2})$ are as in corollary~\ref{ckc.10},
then let $[ \Id +S]_{\infty}$ denote the homotopy class of 
\eqref{ckc.15} in
$[ V\; ; \; \GL(\infty,\bbC)]$.
\label{ckc.18}\end{definition}
The homotopy class $[\Id+ S]_{\infty}$ depends a priori on three choices:
\begin{enumerate}
\item  The choice of an explicit trivialization of $F$, 

\item The choice of a spectral section as in corollary~\ref{ckc.10},

\item The choice of a family of Laplacians $\Delta_{\pa X /Y}$.

\end{enumerate}

In the next three lemmas, we will show that in fact it is independent of 
these three choices.

\begin{lemma}
For $\Delta_{\pa X/Y}$ and $(\Pi,R_{1},R_{2})$ fixed,
the homotopy class $[\Id + S]_{\infty}$ does not depend on the way $F$ is 
trivialized. 
\label{ckc.19}\end{lemma}
\begin{proof}
Suppose $(\Id+S_{\Pi}):V\to \GL(\infty,\bbC)$ and 
$(\Id+S_{\Pi}'):V\to \GL(\infty,\bbC)$ arise from two different
trivializations of $F$.  Then,
\[
             \Id+S_{\Pi}'= M(\Id + S_{\Pi})M^{-1}
\]
for some smooth map $M:V\to\GL(\infty,\bbC)$ which is a pull-back
of a map from $Y$ to $\GL(\infty,\bbC)$.  But at the homotopy level,
the product  $[ a ]\circ [ b ]= [a  b ]$ is commutative (see for instance
 lemma $2.4.6$ in \cite{Atiyah}), so
\[
\begin{aligned}
   \left[ \Id+ S_{\Pi}'\right]  
       &= [ M(\Id+S_{\Pi})M^{-1}]=[ M ] \circ[ \Id+ S_{\Pi} ] \circ [ M^{-1} ]
  \\
                  &= [\Id +S_{\Pi} ] \circ [ M ]\circ [ M^{-1} ]= 
                                    [\Id +S_{\Pi} ] \circ [ \Id ]
  \\
                  &= [ \Id +S_{\Pi} ].
\end{aligned}
\]
\end{proof}

\begin{lemma}
Assume the family of Laplacians $\Delta_{\pa X /Y}$ is fixed and let
$(\Id+S)$ be as in lemma~\ref{ckc.8}.  If $(\Pi_{1},R_{1},T_{1})$ and
$(\Pi_{2},R_{2},T_{2})$ are two spectral sections as in lemma~\ref{ckc.8},
then $\Id+S_{\Pi_{1}}$ and $\Id+ S_{\Pi_{2}}$ define the same homotopy class
in $[ V \; ; \; \GL(\infty,\bbC)]$.
\label{ckc.20}\end{lemma}
\begin{proof}
For $k\in\bbN$, set $R'_{k}=\max\{T_{1},T_{2}\}+k$.  
Then proposition~\ref{ss.4}
asserts the existence of a sequence of spectral sections 
$\{(\Pi_{k}',R_{k}',T_{k}')\}_{k\in\bbN}$ for the family of Laplacians 
$\Delta_{\pa X/Y}$ on $Y$.  By construction,
\[
      \Pi_{k}'\Pi_{i}\Pi_{k}'=\Pi_{i}\quad \mbox{for}\;\; i\in\{1,2\},\; 
k\in\bbN,
\]
so the estimate of lemma~\ref{ckc.8} applies as well to all the projections
$\Pi_{k}'$:
\[
\| S_{b}-\Pi_{k,b}'S_{b}\Pi_{k,b}' \| < \frac{1}{2\| (\Id +S_{b})^{-1}\|},
\quad \forall b\in V,\quad\forall k\in\bbN.          
\]
Moreover, since $R_{k}'\underset{k\to\infty}{\longrightarrow} \infty$,  we know
from \eqref{ckc.7}, the fact that $\Id+S$ is asymptotic to the identity and 
the compactness of $Y\times[0,1]$ that in the 
$\CI$-topology,
\[
   \Pi_{k,b}'(\Id + S_{b} +t(\Pi_{i,b} S_{b}\Pi_{i,b}-S_{b}))\Pi_{k,b}'
    +\Id-\Pi_{k,b}'
\underset{k\to\infty}{\longrightarrow}
   \Id+S_{b} +t(\Pi_{i,b}S_{b}\Pi_{i,b}-S_{b})
\]
uniformly in $(b,t)\in V\times [0,1]$, $i\in\{1,2\}$. By the proof
of corollary~\ref{ckc.10}, the limits of these two sequences are homotopies
of sections of $\pi^{*}\G$ asymptotic to the identity. Thus, taking
$k\in\bbN$ large enough, we can assume that 
\[
  t \mapsto\Pi_{k,b}'(\Id + S_{b} +t(\Pi_{i,b} S_{b}\Pi_{i,b}-S_{b}))\Pi_{k,b}'
    +\Id-\Pi_{k,b}',\quad t\in[0,1],  
\]
are homotopies of sections of $\pi^{*}\G$ asymptotic to the identity 
for $i\in\{1,2\}$.   
Then, if $F_{i}$ and $F_{k}'$ denote the range of $\Pi_{i}$ and $\Pi_{k}'$,
we see that 
\begin{equation}
   t\mapsto \Pi_{k}'( \Id+S + t(\Pi_{i}S\Pi_{i}-S))\Pi_{k}',\quad t\in[0,1],
\label{ckc.21}\end{equation}
is a homotopy of sections of $\GL(F_{k}',\bbC)$ between 
$\Pi_{k}'(\Id+S)\Pi_{k}'$
and $\Pi_{k}'(\Id +\Pi_{i}S\Pi_{i})\Pi_{k}'$ for $i\in \{1,2\}$.  It is 
possible that the complements of the $F_{i}$ in $F_{k}'$ are 
not trivial bundles, but
by considering a spectral section $(\Pi_{j},R_{j}',T_{j}'')$ with 
$R_{j}'>T_{k}'$ and $j$ large enough, we can add a trivial bundle  to $F_{k}'$ 
so that the complements of the $F_{i}$ become trivial\footnote{For instance, we
can add
a trivial bundle isomorphic to $F_{1}\oplus F_{2}$.} and so that
$(\Pi_{k},R_{k},T_{j}'')$ is still a spectral section.  So without loss
of generality, we can assume that the complements of the $F_{i}$ in 
$F_{k}'$ are trivial.  Then, we see from lemma~\ref{ckc.19} that
\[
\begin{aligned}
\left[ \Id+ S_{\Pi_{1}} \right] &= [ \Id+ (\Pi_{1}S\Pi_{1})_{\Pi_{k}'} ] \\
          &=  [ \Id + S_{\Pi_{k}'} ] \quad \mbox{using \eqref{ckc.21} with} 
\;i=1, \\
              &= [ \Id+ (\Pi_{2}S\Pi_{2})_{\Pi_{k}'} ] \quad 
                            \mbox{using \eqref{ckc.21} with}\; i=2, \\
              &= [ \Id+ S_{\Pi_{2}} ].
\end{aligned}
\]
\end{proof}
\begin{lemma}
The homotopy class $[\Id+S]_{\infty}$ of definition~\ref{ckc.18} 
does not depend
on the choice of the family of Laplacians $\Delta_{\pa X/Y}$.
\label{ckc.22}\end{lemma}
\begin{proof}
Let $\Delta_{0}$ and $\Delta_{1}$ be two differents families of Laplacians for 
the fibration $\Phi:\pa X\to Y$ and the complex vector bundle $E$, arising
respectively from the metrics $g_{0}^{T(\pa X/Y)}, g_{0}^{E}$ and 
$g_{1}^{T(\pa X/Y)}, g_{1}^{E}$, and connections $\nabla^{E}_{0}$ and
$\nabla^{E}_{1}$.  We need to show that they lead to 
the same homotopy class in $[V \; ; \; \GL(\infty,\bbC) ]$.

By considering the homotopies of metrics
\[
  t\mapsto g_{t}^{T(\pa X/Y)}=(1-t)g_{0}^{T(\pa X/Y)} + tg_{1}^{T(\pa X/Y)},
\quad t\in[0,1],
\]
and
\[
  t\mapsto g_{t}^{E}=(1-t)g_{0}^{E} + tg_{1}^{E},
\quad t\in[0,1],
\]
and the homotopy of connections
\[
  t\mapsto (1-t)\nabla^{E}_{0} + t\nabla^{E}_{1}, \quad t\in[0,1],
\]
one gets an associated homotopy of families of Laplacians
\[
   t\mapsto \Delta_{t},\quad t\in [0,1],
\]
where $\Delta_{t}$ is defined using the metrics $g_{t}^{T(\pa X/Y)}$ and
$g_{t}^{E}$ and the connection $\nabla^{E}_{t}$.  We can think of 
$\Delta_{t}$ as a family of
Laplacians parametrized by $Y\times [0,1]$.  A moment of
reflection reveals that the proofs of lemma~\ref{ckc.8} and 
corollary~\ref{ckc.10} apply equally well
if $Y$ is replaced by 
$Y\times [0,1]_{t}$
and $\Id +S$ is replaced by the section $\Id + p^{*}S$, where 
$p:V\times [0,1]_{t}\to V$ is the 
projection on the left factor.  This leads to a spectral section 
$(\Pi,R_{1},R_{2})$ satisfying the estimate of lemma~\ref{ckc.8} and such 
that 
\[
       t\mapsto \Id + p^{*}S + t(\Pi(p^{*}S)\Pi- p^{*}S), \quad t\in[0,1]
\]
is a homotopy of sections of $p^{*}\pi^{*}\G$ between
$(\Id + p^{*}S)$ and $(\Id + \Pi(p^{*}S)\Pi)$.  If $\Pi_{t}$ denotes
the restriction of $\Pi$ to the slice $Y\times\{t\}$,
then $(\Pi_{0},R_{1},R_{2})$ and $(\Pi_{1},R_{1},R_{2})$ are spectral 
sections for $\Delta_{0}$ and $\Delta_{1}$ which satisfy the conclusions
of lemma~\ref{ckc.8} and corollary~\ref{ckc.10}.  Moreover, using 
lemma~\ref{ckc.19}, we see from the homotopy
\[
   t \mapsto \Id + \Pi_{t}S\Pi_{t}, \quad t\in[0,1],
\]
that $[\Id + S_{\Pi_{0}}]_{\infty}= [\Id + S_{\Pi_{1}}]_{\infty}$.  From 
lemma~\ref{ckc.20},
we conclude that the homotopy class of definition~\ref{ckc.18} is the same 
for $\Delta_{0}$ and $\Delta_{1}$.
\end{proof}

From lemmas \ref{ckc.19}, \ref{ckc.20} and \ref{ckc.22}, we get the following

\begin{proposition}
The homotopy class $[\Id+S]_{\infty}$ of definition~\ref{ckc.18} is 
well-defined, that
is, it is independent of the three choices involved in its construction.
\label{ckc.23}\end{proposition}

\begin{definition}
Let $[V\; ; \; \pi^{*}\G]$ denote the set of homotopy classes of smooth
sections of $\pi^{*}\G$ asymptotic to the identity.  If $\Id+S$  is a section
of $\pi^{*}\G$ asymptotic to the identity, let $[\Id+S]$ denotes its homotopy
class in $[V\; ; \; \pi^{*}\G]$.
\label{ckc.34}\end{definition}
\begin{lemma}
If $\Id+S$ is a section of $\pi^{*}\G$ asymptotic to the identity, then
$[\Id+S]_{\infty}$ in $[V\; ; \;\GL(\infty,\bbC)]$ only depends on the 
homotopy class of $\Id+S$ in $[V\; ; \; \pi^{*}\G]$.
\label{ckc.35}\end{lemma}
\begin{proof}
Assume that $t\mapsto \Id+S_{t}$, $t\in[0,1]$, is a homotopy of sections
of $\pi^{*}\G$ asymptotic to the identity.  We need to show that 
$[\Id+S_{0}]_{\infty}=[\Id+S_{1}]_{\infty}$.  Since $[0,1]$ is compact,
a moment of thought reveals that the proofs of lemma~\ref{ckc.8} and
corollary~\ref{ckc.10} work equally well if one replaces $V$ by the 
bundle $[0,1]\times V\to Y$.  What one gets is that there exists a spectral
section $(\Pi,R_{1},R_{2})$ for $\Delta_{\pa X/Y}$ on $Y$ so that 
\[
      \| S_{t}(b)-\Pi_{b}S_{t}(b)\Pi_{b}\| < 
\frac{1}{2\| (\Id+S_{t}(b))^{-1}\|}, \quad \forall b\in V, t\in [0,1],
\]
which implies that $\Id+S_{t}$ and $\Id+\Pi S_{t}\Pi$  are homotopic
as smooth homotopies of smooth sections of $\pi^{*}\G$ asymptotic 
to the identity.  In particular, after trivializing the range of $\Pi$, 
$\Id+\Pi S_{0}\Pi$ and $\Id +\Pi S_{1}\Pi$ define the same
homotopy class in $[V\; ; \; \GL(\infty,\bbC)]$, and we conclude 
that 
\[
        [\Id+S_{0}]_{\infty}=[\Id+S_{1}]_{\infty}.
\]
\end{proof}
Thus, we see that definition~\ref{ckc.18} gives us a canonical map
\begin{equation}
      \begin{array}{rrcl}
       I_{\infty}: & [V\; ;\; \pi^{*}\G] & \to & [V\; ;\; \GL(\infty,\bbC)] \\
                   & [\Id+S] & \mapsto & [\Id+S]_{\infty} \; .
      \end{array}
\label{ckc.36}\end{equation}
\begin{proposition}
The canonical map $I_{\infty}$ is an isomorphism of sets.
\label{ckc.37}\end{proposition}
\begin{proof}
The proof of the surjectivity of $I_{\infty}$ reduces to  
the existence
of a spectral section with range a vector bundle of arbitrary large rank.
But the existence of such a spectral section is an easy consequence of 
proposition~\ref{ss.4} together with the compactness of $Y$ and the 
continuity of the spectrum as a set-valued function.

For the proof of the injectivity of $I_{\infty}$, suppose 
$[\Id+S_{0}]_{\infty}=[\Id+S_{1}]_{\infty}$.  Let $\Pi_{0}$ and 
$\Pi_{1}$ be spectral sections as in lemma~\ref{ckc.8} that can be used
to define $[\Id+S_{0}]_{\infty}$ and $[\Id+S_{1}]_{\infty}$ respectively.
Let $\Pi$ be another spectral section such that $\Pi\Pi_{i}\Pi=\Pi_{i}$ for
$i\in\{0,1\}$ (cf. proof of lemma~\ref{ckc.20}).  Then $\Pi$ can also
be used to define $[\Id+S_{0}]_{\infty}$ and $[\Id+S_{1}]_{\infty}$.  Taking
$\Pi$ to have a larger rank if necessary, we can assume that there is 
a homotopy
\[
          t\mapsto \Id + \Pi S_{t}\Pi, \quad t\in[0,1],
\]
through sections of $\pi^{*}G^{-\infty}$ asymptotic to the identity.  From
corollary~\ref{ckc.10}, we then deduce that $[\Id+S_{0}]=[\Id+S_{1}]$
in $[V\; ; \; \pi^{*}\G]$.  This shows that the map
$I_{\infty}$ is injective.
\end{proof}

Before discussing the associated $K$-class, let us concentrate on the
second case.

\subsection{The case where $\dim Z=0$}

The case where $\dim Z=0$ turns out to be significantly easier.  This is
somehow the reverse procedure of the case where $\dim Z>0$, that is, instead
of reducing an infinite  bundle to a trivial vector bundle, we will enlarge 
a finite vector bunlde to make it trivial.  The case $\dim Z=0$ evidently
includes scattering operators, but it means more generally that 
the fibration $\Phi:\pa X\to Y$ is a finite covering.

\begin{definition}
Let $\E\to Y$ be the complex vector bundle on $Y$ with fibre 
at $y\in Y$  given by
\[
     \E_{y}=\bigoplus_{z\in\Phi^{-1}(y)} E_{y}.
\]
\label{ckc.24}\end{definition}  

\begin{definition}
Let $\GL(\E,\bbC)$ be the bundle over $Y$ with fibre at $y\in Y$ given by
$\GL(\E_{y},\bbC)$.  Let $\pi^{*}\GL(\E,\bbC)\to \fcn^{*}Y$ be the pull-back
of $\GL(\E,\bbC)$ on $\fcn^{*}Y$.
\label{ckc.25}\end{definition}
An immediate consequence of remark~\ref{fco.30} and proposition~\ref{fco.18}
is the following 

\begin{corollary}
The indicial family $(\Id+\widehat{A}_{0})$ of a Fredholm operator 
$(\Id+A)\in\Frc(X;E)$ is a section of the bundle $\pi^{*}\GL(\mathcal{E};\bbC)$
which is asymptotic to the identity.
\label{ckc.26}\end{corollary}

Let $\F\to Y$ be a complex vector bundle over $Y$ such that $\E\oplus\F$ is
trivial.  Such a bundle always exists, see for instance corollary $1.4.14$
in \cite{Atiyah}.  If $n$ is the rank of $\E\oplus \F$, so that 
$\E\oplus\F\cong \underline{\bbC}^{n}$, then
\begin{equation}
\begin{aligned}
  \pi^{*}\GL(\E,\bbC)&\subset \pi^{*}\GL(\E\oplus\F,\bbC)\\
          &\cong
           \pi^{*}\GL(\underline{\bbC}^{n},\bbC)\\
           &=\GL(n,\bbC)\times\fcn^{*}Y \\
           &\subset \GL(\infty,\bbC)\times\fcn^{*}Y.
\end{aligned}
\label{ckc.27}\end{equation}

\begin{definition}
If $(\Id+S)$ is the indicial family of some operator in $\Frc(X;E)$,
let $(\Id +S_{\F}):\fcn^{*}Y\to \GL(\infty;\bbC)$ be the associated map under
the series of inclusions \eqref{ckc.27} and let $[ \Id +S ]_{\infty}
\in [\fcn^{*}Y\; ;\; \GL(\infty;\bbC) ]$ denotes its homotopy class.
\label{ckc.28}\end{definition}

As in the previous case, we want to show that the homotopy class 
$[ \Id +S ]_{\infty}$ does not depend on the choices made to define it, 
namely the 
choice of an explicit trivialization of $\E\oplus\F$ and the choice of 
$\F$.

\begin{proposition}
The homotopy class $[ \Id+ S]_{\infty}$ of definition~\ref{ckc.28} 
is well-defined,
 that is, it does not depend on the choice of $\F$ and on the way
$\E\oplus \F$ is trivialized.
\label{ckc.31}\end{proposition}
\begin{proof}
First, the homotopy class $[\Id + S ]_{\infty}$  
does not depend
on the way $\E\oplus\F$ is trivialized.
It is the same proof as lemma~\ref{ckc.19}.

Secondly, the homotopy class $[ \Id + S ]_{\infty}$  
does not depend
on the choice of the vector bundle $\F$ such that $\E\oplus\F$ is trivial.

Indeed,
let $\Gv\to Y$ be another complex vector bundle such that $\E\oplus\Gv$ is 
trivial.  Then consider the trivial vector bundle $\Hv= \E\oplus F\oplus\E'
\oplus\Gv$, where $\E'\cong \E$.  The bundles $\E\oplus\F$ and $\E\oplus\Gv$
are subbundles of $\Hv$, and to be more precise, we consider the inclusions
\[
   \begin{array}{rrcl}
    i_{1}: & \E\oplus \F & \hookrightarrow & \E\oplus F\oplus\E'\oplus\Gv \\
           &  (e,f)  & \mapsto             & (e,f,0,0)
   \end{array}, \quad
\begin{array}{rrcl}
    i_{2}: & \E\oplus \Gv & \hookrightarrow & \E\oplus F\oplus\E'\oplus\Gv \\
           &  (e,g)  & \mapsto             & (e,0,0,g)
   \end{array}.
\]
Clearly, the complements of $\E\oplus\F$ and $\E\oplus \Gv$ in $\Hv$ are 
trivial vector bundles, so this means that 
$\F$ and $\F\oplus\E'\oplus\Gv$ lead to the same homotopy class, and similarly
for $\Gv$ and $\F\oplus\E'\oplus\Gv$.  In particular $\F$ and $\Gv$ lead
to same homotopy class.
\end{proof}

In this case, it is straightforward to see that the homotopy class
$[\Id+S]_{\infty}$ only depends on the homotopy class of $\Id+S$
in $[\fcn^{*}Y\; ; \; \pi^{*}\GL(\E,\bbC)]$.  So definition~\ref{ckc.28}
gives also rise to a canonical map
\begin{equation}
\begin{array}{rrcl}
    I_{\infty}: & [\fcn^{*}Y\; ; \; \pi^{*}\GL(\E,\bbC)] & \to & 
  [\fcn^{*}Y\; ; \; \GL(\infty, \bbC)] \\
                & [\Id+S] & \mapsto & [\Id+S]_{\infty}\; .
\end{array}
\label{ckc.39}\end{equation}
However, it is not an isomorphism in general, but all classes arise this way
if arbitrary bundles are admitted.

\subsection{The associated $K$-class}

Now that in both cases $\dim Z>0$ and $\dim Z=0$ we can associate a 
well-defined homotopy class 
\[
        [ \Id+ \widehat{A}_{0} ]_{\infty}\in [ \fcn^{*}Y \; ; \; 
                                                        \GL(\infty;\bbC) ]
\]
to the indicial family $(\Id+ \widehat{A}_{0})$ of a Fredholm operator
$(\Id+A)$ in $\Frc(X;E)$, let us try to reinterpret this homotopy class
in terms of $K$-theory.  We will follow the notation of \cite{Atiyah}.
Recall that the space $\GL(\infty,\bbC)$ is a 
classifying space for odd $K$-theory\footnote{See for instance lemma $2.4.6$ in
\cite{Atiyah}}, so we have the correspondence
\[
   [\fcn^{*}Y\; ;\;\GL(\infty,\bbC) ] \cong \widetilde{K}^{0}(\bbS^{1}\wedge 
Y^{\fcn^{*}Y}),
\]
where $Y^{\fcn^{*}Y}$ is the \textbf{Thom space} of $\fcn^{*}Y$, in other words 
its one point compactification, and 
\[
      \bbS^{1}\wedge Y^{\fcn^{*}Y}= (\bbS^{1}\times Y^{\fcn^{*}Y})/ 
    (\bbS^{1}\times\{\infty\}\cup \{1\}\times Y^{\fcn^{*}Y})
\]
is the \textbf{reduced suspension} of $Y^{\fcn^{*}Y}$, $\bbS^{1}\subset\bbC$ 
being the unit circle.  Since $\fcn^{*}Y\cong T^{*}Y\times\bbR$, we see that
\[
     Y^{\fcn^{*}Y}\cong \bbS^{1}\wedge Y^{T^{*}Y} = 
                                         (\bbS^{1}\times Y^{T^{*}Y})/
                      (\bbS^{1}\times\{\infty\}\cup \{1\}\times Y^{T^{*}Y}),
\]
the factor ``$\bbR$'' being effectively turned into a suspension.  This means
that 
\begin{equation}
\begin{aligned}
\left[ \fcn^{*}Y; \GL(\infty,\bbC)\right] 
&\cong \widetilde{K}^{0}(\bbS^{1}\wedge \bbS^{1}
\wedge Y^{T^{*}Y}) \\
             &= \widetilde{K}^{-2}(Y^{T^{*}Y}) \\
             &\cong \widetilde{K}^{0}(Y^{T^{*}Y}) \quad \mbox{by the  
                                            periodicity theorem} \\
             &=\K(T^{*}Y),
\end{aligned}
\label{ckc.32}\end{equation} 
where $\K(T^{*}Y)$, the \textbf{$K$-theory with compact support} of $T^{*}Y$,
is by definition equal to $\widetilde{K}^{0}(Y^{T^{*}Y})$.
Thus, we see from \eqref{ckc.32} that the homotopy class 
$[\Id + \widehat{A}_{0} ]$ in $[\fcn^{*}Y\; ; \; \GL(\infty,\bbC) ]$ defines
a $K$-class in $\K(T^{*}Y)$.

\begin{definition}
If $(\Id+ \widehat{A}_{0})$ is the indicial family of a Fredholm operator 
$(\Id+A)$ in $\Frc(X;E)$, let $\kappa(\Id+A)\in \K(T^{*}Y)$ denote the 
$K$-class associated to the homotopy class $[\Id+ \widehat{A}_{0} ]_{\infty}$
 under
the identification $[ \fcn^{*}Y; \GL(\infty,\bbC)]\cong \K(T^{*}Y)$ given
by \eqref{ckc.32}.
\label{ckc.33}\end{definition}

\section{The Index in terms of $K$-theory}\label{im.0}

Having associated a $K$-class $\kappa(\Id+A)\in \K(T^{*}Y)$ to a 
Fredholm operator $(\Id + A)$ in $\Frc(X;E)$, the obvious guess is that the 
index of $\Id+A$ is given by
\begin{equation}
               \ind(\Id +A)= \ind_{t}(\kappa(\Id+A))
\label{im.1}\end{equation}
where
\[
                \ind_{t}:\K(T^{*}Y)\to \bbZ
\]
is the topological index introduced by Atiyah and Singer in 
\cite{Atiyah-Singer1}.  This is indeed the case.  We will prove 
\eqref{im.1} by reducing the computation of the index to the case of 
a scattering operator.

To the manifold $X$ with fibred cusp structure given by the fibration
$\Phi:\pa X\to Y$ and some defining function $x$ of $\pa X$, let us associate
the manifold with boundaries
\begin{equation}
    \W = [0,1]_{t}\times Y, \quad \pa\W =Y_{0}\cup Y_{1},
\label{im.2}\end{equation}
where $Y_{0}=\{0\}\times Y\subset \W$ and $Y_{1}=\{1\}\times Y\subset \W$.
The boundary $\pa \W$ has an obvious defining function given by $t$ near
$Y_{0}$ and by $(1-t)$ near $Y_{1}$, where $t\in[0,1]$ is the variable of 
the left factor of $[0,1]_{t}\times Y$.  It will have a fibred cusp 
structure if one considers the trivial fibration given by the identity
\[
                  \Id: \pa\W\to \pa\W.
\]
If $F\to \W$ is a complex vector bundle on $\W$, then we can define the 
algebra of fibred cusp operators $\fci^{*}(\W;F)$ acting on 
sections of $F$.  From the previous section, we know that a Fredholm
operator in $\Fri(\W;F)$ gives rise to a well define element in
\[
             \K(T^{*}(\pa \W))\cong \K(T^{*}Y_{0})\oplus \K(T^{*}Y_{1}). 
\]
\begin{definition}
Let $i_{1}:\K(T^{*}Y)\hookrightarrow \K(T^{*}(\pa\W))$ be the natural 
inclusion given by
\[
\begin{array}{rrcl}
  i_{1}: & \K(T^{*}Y) & 
\hookrightarrow & \K(T^{*}Y_{0})\oplus \K(T^{*}Y_{1}) \\
         &    a  & \mapsto & (0,a).
\end{array}
\]
\label{im.3}\end{definition}
\begin{lemma}
If $(\Id +A)\in \Frc(X;E)$ is such that the indicial family 
$\widehat{A}_{0}$,
seen as a section of $\pi^{*}\mathcal{P}^{-\infty}$, lies in some 
subbundle $\pi^{*}\End(F)\subset \pi^{*}\mathcal{P}^{-\infty}$, where
$F\to Y$ is some trivial finite rank compex vector bundle, then 
 any scattering operator $(\Id+ S)\in\Fri(\W;F)$ such that
\[
                       \kappa(\Id+S)= i_{1}\kappa(\Id+A), \;\mbox{in}\;
   \K(T^{*}\W)
\]
has the same index as $\Id+A$.
\label{im.4}\end{lemma}
\begin{proof}
By proposition~\ref{itdf.7}, we can assume that $\no'(A)=\widehat{A}_{0}$, that
is, $\widehat{A}_{k}=0$ for $k\in\bbN$.  From proposition~\ref{itdf.11}
and formula \eqref{itdf.20}, we see that a parametrix $\Id+B$ of $\Id+A$ 
with full normal operator
\[
        \widehat{\no(B)}=\sum_{k=0}^{\infty}x^{k}\widehat{B}_{k}
\]
is such that $\widehat{B}_{k}$ is a section of $\pi^{*}\End(F)
\subset\pi^{*}\mathcal{P}^{-\infty}$.  According to proposition~\ref{itdf.6},
the index of $\Id+A$ is then given by 
\begin{equation}
    \ind(\Id+A) = \frac{1}{(2\pi)^{l+1}}\int_{\fcn^{*}Y}  
\Tr_{\Phi^{-1}(y)}
(\widehat{(D_{\log x}A B)}_{l+1}(y,\tau,\eta)) dy d\eta d\tau \; .
\label{im.13}\end{equation}
But from lemma~\ref{itdf.12},
\[
\widehat{\no'(D_{\log x}A)}=
\sum_{k=1}^{\infty} \frac{x^{k}D_{\tau}^{k}\widehat{A}_{0}}{k}.
\]
Thus, since $\widehat{B}_{k}, D_{\tau}^{k}\widehat{A}_{0}$ are sections
of $\pi^{*}\End(F)\subset\pi^{*}\mathcal{P}^{-\infty}$ for all 
$k\in\bbN_{0}$, this means the trace $\Tr_{\Phi^{-1}(y)}$ in \eqref{im.13}
can be replaced by the trace on $F_{y}$,
\begin{equation}
    \ind(\Id+A) = \frac{1}{(2\pi)^{l+1}}\int_{\fcn^{*}Y}  
\Tr_{F_{y}}
(\widehat{(D_{\log x}A B)}_{l+1}(y,\tau,\eta)) dy d\eta d\tau \; .
\label{im.14}\end{equation}
But applying proposition~\ref{itdf.6} to scattering operators, we see from 
proposition~\ref{itdf.11} that \eqref{im.14} also gives the index of an 
operator $(\Id+S)\in\Fri(\W;F)$ such that $\widehat{\noi'(S)}=0$ on
$\fcis^{*}Y_{0}$ and $\widehat{\noi'(S)}=\widehat{A}_{0}$ on
$\fcis^{*}Y_{1}$.  Clearly, $\kappa(\Id+S)=i_{1}\kappa(\Id+A)$.  More 
generally, from formula \eqref{im.14}, definition~\ref{ckc.28} and 
proposition~\ref{itdf.6}, any $(\Id+Q)\in\Fri(\W;F)$ 
such that $\kappa(\Id+Q)=i_{1}\kappa (\Id+A)$ has the same index as
$\Id+S$ and $\Id+A$.
\end{proof}

\begin{lemma}
Let $(\Id+ A)\in \Frc(X;E)$ be a Fredholm operator and let $\kappa(\Id+A)$
be its associated $K$-class in $\K(T^{*}Y)$.  Then, for $F$ a trivial vector 
bundle with sufficiently large rank, there exists a scattering operator
$\Id+S\in \Fri(\W;F)$ with associated $K$-class given by 
$i_{1}\kappa(\Id+A)\in\K(T^{*}(\pa \W))$ and such that 
$\ind(\Id+A)=\ind(\Id+S)$.   
\label{im.5}\end{lemma}
\begin{proof}
The proof relies on lemma~\ref{im.4}, but is slightly different depending
on whether or not $\dim Z>0$.

First assume that $\dim Z>0$.  Then let $(\Pi,R_{1},R_{2})$ be a spectral 
section as in lemma~\ref{ckc.8} for the indicial family $\Id+\widehat{A}_{0}$.
By corollary~\ref{ckc.10}, $(\Id+\Pi\widehat{A}_{0}\Pi)$ is homotopic
to $(\Id+\widehat{A}_{0})$.  Let $P\in\fc^{-\infty}(X;E)$ be such that 
$\widehat{P}_{0}=\Pi\widehat{A}_{0}\Pi$.  Then by proposition~\ref{fco.18}
and proposition~\ref{itdf.7}, $(\Id+P)\in\Frc(X;E)$ has the same
index as $\Id+A$ and it defines the same $K$-class in $\K(T^{*}Y)$.  By taking
$F\to Y$ to be the range of $\Pi$ on $Y$, we can then apply 
lemma~\ref{im.4} to $\Id+P$.  Thus there exists 
$(\Id+S)\in\Fri(\W;F)$ such that $\ind(\Id+S)=\ind(\Id+P)$ and 
$\kappa(\Id+S)=i_{1}\kappa(\Id+P)$.  So $\Id+S$ is the desired operator.

If $\dim Z=0$ let $\E\to Y$ be as in definition~\ref{ckc.24}.  From
lemma~\ref{im.4}, there exists $(\Id+S)\in\Fri(\W;\E)$ having 
the same index as $\Id+A$ and having $i_{1}\kappa(\Id+A)$ as an associated
$K$-class.  If $\E$ is not trivial, let $\Gv\to \W$ be another complex 
vector bundle such that $F=\E\oplus\Gv$ is trivial and let 
$\Id+S'\in\Fri(\W;F)$ be the operator which acts as $(\Id+ S)$ on sections
of $\E$ and as the identity on sections of $\Gv$.  Clearly,  $\Id+S'$ has
the same index as $\Id+S$ and they define the same $K$-class.  So
$\Id+S'$ is the desired operator. 
   
\end{proof}

This reduces the problem to the computation of the index of a scattering
operator.  In section $6.5$ of \cite{MelroseGST}, a general topological 
formula for the index of scattering operators is derived.  For the convenience
of the reader, we will adapt the discussion that can be found there to our
particular context.  

Given $(\Id+S)\in\Fri(\W;F)$ as in lemma~\ref{im.5} with 
$F=\underline{\bbC}^{n}$, the indicial family $\Id+\widehat{S}_{0}$ of 
$\Id+S$ can be seen as a map
\begin{equation}
       \Id+\widehat{S}_{0}: \fcis^{*}\pa\W \to \GL(n,\bbC).    
\label{im.6}\end{equation}
The homotopy class
of this map gives rise to the associated $K$-class in $\K(T^{*}\pa\W)$ under 
the identification \eqref{ckc.32}.  Let 
$\Vy$ be the double of $\W$ obtained by taking two copies
$\W$ and $\W'$ of $\W$ and identifying $Y_{i}$ of the first copy with 
$Y_{i}'$ of the second copy for $i\in\{0,1\}$.  Thus, 
$\Vy\cong \bbS^{1}\times Y$.  

Inside $\Vy$, $\fcis^{*}\pa\W$ naturally 
identifies with the conormal bundle $N^{*}(T^{*}\pa\W)$ of $T^{*}\pa\W$.
Via the clutching construction (see p.75 in \cite{Atiyah}), the map   
\eqref{im.6} defines a $K$-class with compact support in $\K(T^{*}\Vy)$.
More generally,

\begin{definition}
For $(\Id+S)\in \Fri(\W;\underline{\bbC}^{n})$, let $c(\Id+S)\in\K(T^{*}\Vy)$
be the element defined by applying the clutching construction to the map
\eqref{im.6}.  
\label{im.7}\end{definition}  
The following proposition is a particular case of theorem $6.4$ in 
\cite{MelroseGST}.

\begin{proposition}
The index of $(\Id+S)\in\Fri(\W;\underline{\bbC}^{n})$ is given by
\[
       \ind(\Id+S)= \ind_{t}(c(\Id+S))
\]
where $\ind_{t}:\K(T^{*}\Vy)\to \bbZ$ is the topological index map.
\label{im.8}\end{proposition}

\begin{proof}
Let us enlarge our space of operators and see $\Id+S$ as a Fredholm 
operator in $\fci^{0}(\W;\underline{\bbC}^{n})$.  Then proposition~\ref{fco.18}
can be reinterpreted as saying that the $K$-class defined by the symbol
$\sigma_{0}(\Id+S)$ is null when restricted to $\left. T^{*}\W\right|_{\pa\W}$.
This means one can smoothly deform $\Id+S$ through Fredholm operators in
$\fci^{0}(\W;\underline{\bbC}^{n})$ to a Fredholm operator 
$P\in \fci^{0}(\W;\underline{\bbC}^{n})$  that acts by multiplication of a 
matrix near $\pa\W$.  This new operator is no longer in 
$\Fri(\W;\underline{\bbC}^{n})$ in general, but it has the same index as 
$\Id+S$.  
A general Fredholm operator $Q\in \fci^{0}(\W;\underline{\bbC}^{n})$ 
also defines
an elment of $\K(T^{*}\Vy)$ but one needs also the symbol.  If 
$\overline{T^{*}\W}$ is the radial compactification of $T^{*}\W$, then the 
indicial family $\widehat{Q}_{0}$ together with the symbol $\sigma_{0}$ form a 
continuous map  
\[
      f:\pa(\overline{T^{*}\W}) \to GL(n,\bbC). 
\]
This map can be be used to defined a relative $K$-class on the double 
$D_{\overline{T^{*}\W}}$ 
of $\overline{T^{*}\W}$ , obtained by identifying two copies of 
$\overline{T^{*}\W}$ at their boundaries.  This class will in fact have 
support inside $T^{*}\Vy\subset D_{\overline{T^{*}\W}}$ .  
Thus this relative $K$-class is really an element of $\K(T^{*}\Vy)$.  When
$Q\in \Fri(\W;\underline{\bbC}^{n})$, it is not hard to see that this 
construction reduces to definition~\ref{im.7}.  

Coming back to $P$, this means that $P$ defines the same element in
$\K(T^{*}\Vy)$ as $c(\Id+S)$.  Now, since $P$ acts as a matrix near
the boundary, it defines  a vector bundle $\E\to \Vy$ via the 
clutching construction and $P$ can be extended to act as the identity
on $\left. \E\right|_{\W'}$.  This defines a pseudodifferential operator
$\widetilde{P}\in \Psi^{0}(\Vy;\E)$ having the same index as $P$.  Moreover,
it is not hard to see that the $K$-class in $\K(T^{*}\Vy)$ defined by the 
symbol of $\widetilde{P}$ is simply $c(\Id+S)$.  By the Atiyah-Singer index
theorem, the index of $\widetilde{P}$ is then given by $\ind_{t}(c(\Id+S))$, 
which completes the proof.

\end{proof}

The formula for the index is now easy to get.  

\begin{theorem}
The index of a Fredholm operator $(\Id+A)\in\Frc(X;E)$ is given by 
\[
             \ind(\Id+A)= \ind_{t}(\kappa(\Id+A))
\]
where $\kappa(\Id+A)\in\K(T^{*}Y)$ is the associated $K$-class of 
definition~\ref{ckc.33} and 
\[\ind_{t}:\K(T^{*}Y)\to\bbZ \]
 is the topological index map of Atiyah and Singer.
\label{im.9}\end{theorem}
\begin{proof}
By lemma~\ref{im.5}, for $n\in\bbN$ large enough, there exists a scattering
operator $(\Id+S)\in\Fri(\W;\underline{\bbC}^{n})$ which has the same
index as $\Id+A$ and such that $\kappa(\Id+S)=i_{1}\kappa(\Id+A)$.   
By the previous proposition, the index of $\Id+S$ is given
by
\[
             \ind(\Id+S)=\ind_{t}(c(\Id+S)),
\]
where $c(\Id+S)\in\K(T^{*}\Vy)$ is defined using the clutching construction.
At $Y_{0}\subset\pa\W$, the indicial family of $(\Id+S)$ is trivial.  
This means the clutching construction only depends on the part of the indicial 
family of $\Id+S$ defined over $\fcis^{*}Y_{1}$. 
Let 
\begin{equation}
  c: [\fcn^{*}Y\; ;\;\GL(n;\bbC)]\to \K(T^{*}\Vy)
\label{im.10}\end{equation}
denotes the map obtained by applying the clutching construction, where $Y$
is identified with $Y_{1}\subset\Vy$.  The
identification
\[
   [\fcn^{*}Y\; ;\;\GL(\infty;\bbC)]\cong[\bbR\times T^{*}Y\; ;\;
\GL(\infty;\bbC)]\cong \widetilde{K}^{-2}(Y^{T^{*}Y})
\]
combined with \eqref{im.10} gives a map
\[
          f:\widetilde{K}^{-2}(Y^{T^{*}Y})\to \K(T^{*}\Vy).
\]
Let us decompose the map $f$ into simpler maps.
If we identify the tangent bundle with the cotangent bundle via some metric,
$T^{*}Y_{1}$ can be seen as a submanifold of $T^{*}\Vy$. Then, 
a tubular neighborhood of $T^{*}Y_{1}$ in $T^{*}\Vy$ 
defines an inclusion 
\[
            i: N(T^{*}Y_{1})\hookrightarrow T^{*}\Vy
\]
of the normal bundle $N(T^{*}Y_{1})$ of $T^{*}Y_{1}$.  This
in turn defines a push-forward map in $K$-theory
\[
    i_{*}:  \K(N(T^{*}Y_{1}))\to \K(T^{*}\Vy).
\]
On the other hand, since 
$N(T^{*}Y_{1})\cong\bbR^{2}\times T^{*}Y$,  we see that 
there is an isomorphism
\[
   j: \widetilde{K}^{-2}(Y^{T^{*}Y})\overset{\sim}{\longrightarrow}
\K( N(T^{*}Y_{1})).
\]
It is not hard to see that the map $f$ is given by
\[
         f=i_{*}\circ j: \widetilde{K}^{-2}(Y^{T^{*}Y})\to \K(T^{*}\Vy).
\]
Combined with the periodicity isomorphism $p:\K(T^{*}Y)\overset{\sim}{\to}
\widetilde{K}^{-2}(Y^{T^{*}Y})$, we see that 
\[
                c(\Id+S)= i_{*}\circ j\circ p(\kappa(\Id+A)).
\]
Now, the composition of maps $j\circ p$ is just the Thom isomorphism
\[
        \varphi: \K(T^{*}Y)\to \K(N(T^{*}Y_{1}))   
\]
applied to the trivial bundle $\underline{\bbC}$ over $T^{*}Y$,
since $N(T^{*}Y_{1})\cong \bbC\times T^{*}Y$.  Thus, one
has 
\[
                 c(\Id+S) = i_{!}\kappa(\Id+A)
\]
where $i_{!}=i_{*}\varphi$.  From the commutativity of the 
diagram\footnote{see p.501 in \cite{Atiyah-Singer1}} 
\[
  \begin{CD}
           \K(T^{*}Y)  @>{i_{!}}>> \K( T^{*}\Vy) \\
      @VV{\ind_{t}}V       @VV{\ind_{t}}V \\
         \bbZ @=   \bbZ
  \end{CD},         
\]
we conclude that
\[
\begin{aligned}
 \ind(\Id+A) &= \ind(\Id+S) \\
             &= \ind_{t}c(\Id+S)= \ind_{t} i_{!}\kappa(\Id+A)\\
             &= \ind_{t}\kappa(\Id+A).
\end{aligned}
\]

\end{proof}

\begin{corollary}
When $\dim Z> 0$, the index map $\ind:\Frc(X;E)\to \bbZ$ is surjective.
\label{im.11}\end{corollary}
\begin{proof}
From theorem~\ref{im.9}, it suffices to prove that the map
\[
       \kappa:\Frc(X;E)\to \K(T^{*}Y)
\]
of definition~\ref{ckc.33} is surjective, which is equivalent to showing that
any homotopy class in $[\fcn^{*}Y\; ;\; \GL(\infty,\bbC)]$ can be realized
through definition~\ref{ckc.18} as $[\Id+A]_{\infty}$ 
for some $(\Id+A)\in\Frc(X;E)$.

So let  $[f]\in[\fcn^{*}Y\; ;\; \GL(\infty,\bbC)]$ be arbritrary.  Then, there
exists $n\in\bbN$ depending on $[f]$ so that the homotopy class $[f]$ has
a smooth representative 
\[
          f:\fcn^{*}Y\to \GL(n,\bbC)\subset \GL(\infty,\bbC).
\]   
Considering some family of Laplacians as in the beginning of 
section~\ref{ckc.0} and using proposition~\ref{ss.4}, there exists
an associated spectral section $(\Pi,R_{1},R_{2})$ with range a trivial
complex vector bundle $\underline{\bbC}^{m}\to \fcn^{*}Y$ with $m\ge n$.  By
taking $n$ larger if needed, we can assume $m=n$.
Thus taking $(\Id+A)\in\Frc(X;E)$ so that 
\[
            \Pi_{b}(\Id+\widehat{A}_{0}(b))\Pi_{b}= f(b), 
                             \quad \forall b\in\fcn^{*}Y,
\] 
\[
    (\Id-\Pi_{b})(\widehat{A}_{0}(b))(\Id-\Pi_{b})=0, \quad \forall b\in 
\fcn^{*}Y,
\]
we see that $[\Id+A]_{\infty}=[f]$.  
\end{proof}

\begin{corollary}
When $\dim Z=0$ the index map $\ind:\Frc(X;\underline{\bbC}^{n})\to \bbZ$
is surjective provided $n\in\bbN$ is large enough.
\label{im.12}\end{corollary}
\begin{proof}
As in the previous corollary, it suffices to show that any homotopy class
in $[\fcn^{*}Y\; ;\; \GL(\infty,\bbC)]$ can be written as $[\Id+A]$ for some
$\Id+A\in \Frc(X;\underline{\bbC}^{n})$.  But  
there exists $n\in\bbN$ depending on the dimension of $Y$ so that 
any homotopy class in $[\fcn^{*}Y\; ;\; \GL(\infty,\bbC)]$ can be represented
by a map from $\fcn^{*}Y$ to $\GL(n,\bbC)$, and clearly such maps can always be
chosen to be the indicial family of some 
Fredholm operator in $\Frc(X;\underline{\bbC}^{n})$.
  
\end{proof}

\section{The Homotopy Groups of $\Gc$}\label{hg.0}

When homotopy groups come into play, the cases where $\dim Z>0$ and 
where $\dim Z=0$ are radically different.  In some sense, the case where
$\dim Z>0$ is to the case where $\dim Z=0$ what stable homotopy groups 
of $\GL(n, \bbC)$ are to the actual homotopy groups of $\GL(n, \bbC)$.  That
is to say, the homotopy groups are much easier to describe in the first case.
What explains this difference is that the bundle $\pi^{*}\G$ is infinite 
dimensional in the case where $\dim Z>0$.  This gives more freedom in the way
one can perform homotopies, which has the effect of reducing the set 
of homotopy classes, and thus of simplifying the picture.

Still, it is possbible to say something in the case where $\dim Z=0$ if one
consents to allow some stabilization.  This goes as follows.  Let $X$ be 
a compact manifold with boundary $\pa X$ together with a fibred cusp structure
given by a finite covering $\Phi:\pa X\to Y$.  Then consider the manifold with
boundary $\bbS^{1}\times X$ with fibred cusp structure given by
\[ 
          \begin{array}{rrcl}
          \Phi_{\bbS^{1}}: & \bbS^{1}\times X & \to & Y \\
                           & (s,p) & \mapsto & \Phi(p)\; .
          \end{array}
\] 
Let $E\to X$ be any complex vector bundle over $X$.  Then $E$ can be seen as
a subbundle of the (trivial) Hilbert bundle $\Ld(S^{1})$ over $X$, and 
Fredholm operators in $\Frc(X;E)$ can be extended by letting them act
as the identity on the complement of $E$ in $\Ld(\bbS^{1})\to X$.  This gives
an embedding
\[
         \Frc(X;E) \subset \mathcal{F}^{-\infty}_{\Phi_{\bbS^{1}}}(\bbS^{1}
\times X; \underline{\bbC}).
\]
In this way, the space $\mathcal{F}^{-\infty}_{\Phi_{\bbS^{1}}}(\bbS^{1}
\times X; \underline{\bbC})$ is a stabilization of $\Frc(X;E)$ which brings
us back to the case where $\dim Z>0$.  To the reader not completely satisfied
by this compromise, let us point out where our argument breaks down when
$\dim Z= 0$.  It is in the failure of being able to prove lemma~\ref{hg.6}, 
which is
an important ingredient in the proof of the surjectivity of the boundary 
homomorphism occurring in the long exact sequence \eqref{hg.21} .

From now on, we will assume that $\dim Z>0$.
\begin{definition}
Let $\Gc(X;E)\subset \Frc(X;E)$ be the group of invertible operators
in $\Frc(X;E)$, that is,
\[
 \Gc(X;E)=\{ \Id +A \;\; |\;\;A\in\fc^{-\infty}(X;E)\;\mbox{and $(\Id+A)$ is
invertible}\}.
\]
\label{hg.7}\end{definition}
We want to compute the homotopy groups of $\Gc(X;E)$.  Our approach will be
similar to the one in \cite{fipomb} in the sense that we will compute
the homotopy groups of $\Gc(X;E)$ out of a long exact sequence of homotopy
groups.  To this end, let us introduce the other spaces involved.

\begin{definition}
Let $\Gd(X;E)$ be the group of operators given by
\[
  \Gd(X;E)= \{ \Id+Q \; | \; Q\in x^{\infty}\fc^{-\infty}(X;E), \; \Id + Q\;
  \mbox{is invertible} \}.
\]
\label{hg.14}\end{definition}
\begin{remark}
By considering the eigensections vanishing in Taylor series at the boundary of 
a Laplacian arising from a fibred
cusp metric (see ($8.1$) in \cite{mazzeo-melrose4}), one can identify
$\Gd(X;E)$ with $\mathcal{G}^{-\infty}$, the group of semi-infinite invertible
matrices satisfying \eqref{ckc.6}.  This means $\Gd(X;E)$ is a classifying
space for odd $K$-theory.  Notice also that $\Gd(X;E)$ does not depend
on the choice of a fibred cusp structure.
\label{hg.15}\end{remark}
\begin{definition}
Let $\Gs(\pa X;E)$ denote the group
\[
  \{ \Id + Q \;\; | \; \; Q\in \fs^{-\infty}(\pa X; E) \; \mbox{ and 
$\Id+Q$ is invertible} \},
\]
that is, the group of invertible $\fcn Y$-suspended smoothing perturbations
of the identity.  Taking the Fourier transform in $\fcn Y$, this group 
naturally identitfies with the group $\Gamma_{\Id}(\fcn^{*}Y; \pi^{*}\G)$ of 
smooth sections of $\pi^{*}\G$ asymptotic to the identity.  This
is the point of view we will adopt.
\label{hg.1}\end{definition} 
\begin{lemma}
The homotopy groups of $\Gs(\pa X;E)$ are:
\[
               \pi_{k}(\Gs(\pa X;E))\cong \left\{ 
\begin{array}{lcl}
   \K(T^{*}Y) & ,& k \,\mbox{even}, \\
    \widetilde{K}^{-1}(Y^{T^{*}Y}) & , & k \, \mbox{odd}.
\end{array} \right.
\]
\label{hg.2}\end{lemma}
\begin{proof}
By proposition~\ref{ckc.37}, we have
\[
\begin{aligned}
\pi_{k}(\Gs(\pa X;E)) &= \pi_{k}(\Gamma_{\Id}(\fcn^{*}Y; \pi^{*}\G)) \\
                      &\cong [ \fcn^{*}Y\times \bbR^{k}\; ; \; \pi^{*}\G ] \\
                      &\cong [T^{*}Y\times \bbR^{k+1}\; ; \; \GL(\infty,\bbC) 
                                                                    ] \\
                      &\cong \widetilde{K}^{-1}(\bbS^{k+1}\wedge Y^{T^{*}Y})
                          =\widetilde{K}^{-k-2}(Y^{T^{*}Y}).  
\end{aligned}
\]
The result then follows from the periodicity theorem of $K$-theory.
\end{proof}
\begin{proposition}
The homotopy groups of $\Frc(X;E)$ are the same as those of $\Gs(\pa X;E)$, 
namely
\[
               \pi_{k}(\Frc(X;E))\cong \left\{ 
\begin{array}{lcl}
   \K(T^{*}Y) & ,& k \,\mbox{even}, \\
    \widetilde{K}^{-1}(Y^{T^{*}Y}) & , & k \, \mbox{odd}.
\end{array} \right.
\]
\label{hg.3}\end{proposition}
\begin{proof}
Consider the short exact sequence
\[
   \Id+ x\fc^{-\infty}(X;E) \hookrightarrow \Frc(X;E)
   \overset{\no}{\longrightarrow} \Gs(\pa X;E).
\]
This is a Serre fibration.  One can see this using Seeley extension
for manifolds with corners and thinking in terms of Schwartz kernels.
The space $\Id+ x\fc^{-\infty}(X;E)$ is obviously contractible.  So we deduce 
from the associated long exact sequence of homotopy groups
that for $k\in \bbN$,
\[
             \pi_{k}(\Frc(X;E))\cong \pi_{k}(\Gs(\pa X; E)).
\]
That $\pi_{0}(\Frc(X;E))\cong \pi_{0}(\Gs(\pa X;E))$ is clear from Seeley
extension for manifolds with corners and the contractibility of 
 $(\Id+ x\fc^{-\infty}(X;E))$.
\end{proof}

\begin{definition}
Let $\Frco(X;E)$ denote the subspace of Fredholm operators in
$\Frc(X;E)$ having  index zero.  Let $\Gso(\pa X;E)$ be 
the subgroup of $\Gs(\pa X;E)$  given by
\[
    \Gso(\pa X;E)=\widehat{\no(\Frco(X;E))},
\]
and let $\Gso(X;E)[[x]]$ denotes the group of operators 
\[
   \Gso(X;E)[[x]]=\widehat{\no'(\Frco(X;E))}   
\]
where the composition law is given by the $*$-product of 
definition~\ref{itdf.10}.
\label{hg.17}\end{definition}
\begin{remark}
Notice that $\Gso(\pa X;E)$ has the same homotopy groups as 
the space $\Gs(\pa X;E)$ but fewer connected components, more precisely,
\[
\pi_{0}(\Gso(\pa X;E))\cong \pi_{0}(\Frco(X;E))\cong \ker\left[\ind_{t}:
\K(T^{*}Y)\to \bbZ\right].
\]  
Moreover, keeping only the
term of order $x^{0}$ gives a deformation retraction of $\Gso(X;E)[[x]]$ onto
$ \Gso(X;E)$, so $\Gso(X;E)[[x]]$ as the same homotopy groups and the same
set of connected components as $ \Gso(X;E)$.
\label{hg.22}\end{remark}
\begin{lemma}
The full normal operator $\no'$ gives rise to a short exact sequence
\begin{equation}
 \Gd(X;E)\hookrightarrow \Gc(X;E) \overset{\no'}{\longrightarrow}
       \Gso(X;E)[[x]]
\label{hg.19}\end{equation}
which is a Serre fibration.
\label{hg.18}\end{lemma}  
\begin{proof}
The injectivity and the exactness in the middle are clear.  For the 
surjectivity, let $(\Id+Q)\in \Gso(X;E)[[x]]$ be arbitrary, and let 
$(\Id+P)\in \Frco(X;E)$ be such that $\widehat{\no'(P)}=Q$.  
By proposition~\ref{fco.32}, the kernel and the cokernel of $\Id+P$ are
contained in $\CId(X;E)$.  Since $\Id+P$ has index zero, there exists
a linear isomorphism $\varphi: \ker P \to \coker P$.  By extending
it to act as $0$ on the orthogonal complement of $\ker P$ in 
$\Ld(X;E)$, it can be seen as an element of $x^{\infty}\fc^{-\infty}(X;E)$,
so that 
\[
   \Id + P +\varphi \in \Gc(X;E), \;\; \widehat{\no'(P+\varphi)}=Q.
\]  
This shows that $\no'$ is surjective.  The proof that \eqref{hg.19} is 
a Serre fibration is basically the same as the proof of lemma $3.5$ in
\cite{fipomb}.
\end{proof}
A Serre fibration has an associated long exact sequence of homotopy groups, so 
we have the following
\begin{corollary}
The Serre fibration \eqref{hg.19} has an associated long exact sequence
of homotopy groups
\begin{multline}
\cdots \overset{\pa}{\longrightarrow} \pi_{k}(\Gd)\longrightarrow \pi_{k}(\Gc)
\longrightarrow \pi_{k}(\Gso)\overset{\pa}{\longrightarrow}
 \pi_{k-1}(\Gd)\longrightarrow\cdots \\
  \cdots\overset{\pa}{\longrightarrow}\pi_{0}(\Gd)\longrightarrow 
\pi_{0}(\Gc)\longrightarrow \pi_{0}(\Gso),
\label{hg.21}\end{multline}
where $\pa $ is the boundary homomorphism and $\Gd$, $\Gc$ and $\Gso$ denote
respectively $\Gd(X;E)$, $\Gc(X;E)$ and $\Gso(\pa X;E)[[x]]$.
\label{hg.20}\end{corollary}

Since we know the homotopy groups of $\Gd(X;E)$ and $\Gso(X;E)$, we will
be able to compute the homotopy groups of $\Gc(X;E)$ provided we 
identify the boundary homomorphism $\pa$.  In fact, it is only necessary
to know that $\pa$ is surjective.  But even with a complete knowledge of 
$\pa$, there would still be an ambiguity concerning the connected 
components of $\Gc(X;E)$.  Let us get rid of this ambiguity immediately.

\begin{lemma}
The set of connected components of $\Gc(X;E)$ is given by
\[
  \pi_{0}(\Gc(X;E))\cong \ker\left[\ind_{t}:\K(T^{*}Y)\to \bbZ\right].
\]
\label{hg.23}\end{lemma}  
\begin{proof}
From proposition~\ref{hg.3} and theorem~\ref{im.9}, it suffices to 
prove that 
\begin{equation}
    \pi_{0}(\Gc(X;E))\longrightarrow \pi_{0}(\Frc(X;E))
\overset{\ind}\longrightarrow \bbZ
\label{hg.24}\end{equation}
is a short exact sequence.  The surjectivity follows from 
corollary~\ref{im.11}.
The exactness in the middle follows using proposition~\ref{fco.32} (cf. the 
proof of the surjectivity in lemma~\ref{hg.18}).  Finally, the injectivity
follows from the long exact sequence \eqref{hg.21}, the fact that 
$\pi_{0}(\Gd(X;E))=\{0\}$ and proposition~\ref{hg.3}.
\end{proof}

To prove that the boundary homomorphism is surjective, we will see 
it as a generalization of the index map $\ind: \pi_{0}(\Frc(X;E))\to\bbZ$ to 
all homotopy groups of $\Frc(X;E)$.  Let us recall from the appendix in
\cite{Atiyah} how can one performs such a generalization.  

\begin{definition}
If $M$ is a compact manifold with basepoint $m_{0}$, we denote 
by $[M\; ; \; \Frc(X;E) ]$  the set
of homotopy classes of continuous maps from $M$ to $\Frc(X;E)$ which send
the basepoint $m_{0}$ to the identity.  
\label{hg.8}\end{definition}
We will define an 
index map
\begin{equation}
\ind : [M\; ; \;\Frc(X;E)]\to K^{0}(M),
\label{hg.9}\end{equation}
where $\widetilde{K}^{0}(M)$ is the reduced even $K$-theory of $M$.  
Let $H$ denote the Hilbert space $\Ld(X;E)$ defined using
some smooth positive density on $X$.  Given a 
continuous map
\[
              \begin{array}{rrcl}
                  T: & M & \to & \Frc(X;E) \\
                     & m & \mapsto & T_{m},
              \end{array}
\]
with $T_{m_{0}}=\Id$, one can show\footnote{see proposition $A5$ 
in \cite{Atiyah}} the existence of a vector space 
$V\subset H$ which is closed
of finite codimension such that
$V\cap \ker T_{m}=\{0\}$  for all $m\in M$, and such that
\[
\bigcup_{m\in M} H/T_{m}(V),
\]
topologized as a quotient space of $M\times H$, is a vector bundle over $M$.
Let us denote this vector bundle by $H/T(V)$.  

\begin{definition}
Let $T:M\to \Frc(X;E)$ and $V\subset H$ be as above.  Then the \textbf{index}
of $T$ is defined to be
\[
               \ind(T)= [ H/V] -[H/T(V)] \in \widetilde{K}^{0}(M).
\]
\label{hg.10}\end{definition}
One can check that the index of $T$ only depends on the homotopy class of 
$T$ and that it does not depend on the choice of $V\subset H$.  Moreover,
if $S:M\to \Frc(X;E)$ is another continuous map, then
\begin{equation}
        \ind(TS)=\ind(T) + \ind(S).
\label{hg.11}\end{equation}
For later convenience, we will give some precision about the way one can 
choose the vector space $V$ used in definition~\ref{hg.10}.  First, notice
that the orthogonal complement $V^{\perp}$ of $V$ in $H$ is naturally
isomorphic to $H/V$.  Moreover, the vector bundle $T(V)^{\perp}\to M$ with
fibre at $m\in M$ given by the orthogonal complement $T_{m}(V)^{\perp}$ of 
$T_{m}(V)$ is naturally isomorphic to the vector bundle $H/T(V)$.  Hence,
the index can also be written
\[
          \ind(T)= [V^{\perp}]-[ T(V)^{\perp}].
\] 
\begin{lemma}
It is possible to choose the vector space $V$ occuring in 
definition~\ref{hg.10} such that 
\[
V^{\perp}\subset \CId(X;E),\; T_{m}(V)^{\perp}\subset\CId(X;E)\;\forall m\in M.
\] 
Moreover, if $\ind(T)=0$, we can also choose $V$ so that the vector
bundles $V^{\perp}$ and $T(V)^{\perp}$ are isomorphic.
\label{hg.25}\end{lemma}
\begin{proof}
According to the proof of proposition $A5$ in \cite{Atiyah}, one can always 
choose $V$ to be of the form
\[
            V= \bigcap_{i=1}^{n} (\ker T_{m_{i}})^{\perp}\; \Rightarrow
\; V^{\perp}=\sum_{i=1}^{n} \ker T_{m_{i}},
\]
for some finite number of points $m_{1},\ldots,m_{n}\in M$.  Thus,
by propositon~\ref{fco.32}, $V^{\perp}\subset \CId(X;E)$.  Moreover,
for all $m\in M$,
\begin{equation}
      v\in T_{m}(V)^{\perp}\; \Rightarrow \; v\in \ker T^{*}_{m}
+T_{m}(V^{\perp}),
\label{hg.26}\end{equation}
so by proposition~\ref{fco.32}, this means $T_{m}(V)^{\perp}\subset
\CId(X;E)$.  

If $\ind(T)=0$, then there exists a trivial bundle $P$ over $M$ so that
\[
                  H/V\oplus P\cong H/T(V)\oplus P.
\]
By taking a subspace $W\subset V$ such that $\dim (V/W)=\rank(P)$, we
have
\[
              H/W \cong H/T(W).
\]
Without loss of generality, we can assume as well that 
$W^{\perp}\subset\CId(X;E)$, which implies as in \eqref{hg.26} that 
$T_{m}(W)^{\perp}\subset\CId(X;E)$ for all $m\in M$. 
\end{proof}

We know from proposition~\ref{hg.3} that the normal operator $\no$ induces
an isomorphism 
\begin{equation}
 \nos: \pi_{k}(\Frc(X;E))\to 
\pi_{k}(\Gs(X;E))=\pi_{k}(\Gso(X;E)[[x]]), \; k\in \bbN.
\label{hg.27}\end{equation}
Moreover, since $\Gd(X;E)$ is a classifying space for odd $K$-theory, we get
an isomorphism via the clutching construction
\begin{equation}
c: \pi_{k}(\Gd(X;E))\to \widetilde{K}^{-2}(\bbS^{k-1})=\widetilde{K}^{0}(
\bbS^{k+1}), \; k\in \bbN.
\label{hg.28}\end{equation}
\begin{proposition}
For $k\in \bbN$, the homomorphism 
\[
c\circ \pa \circ \nos: \pi_{k}(\Frc(X;E))
\to \widetilde{K}^{0}(\bbS^{k})
\]
is the index map of definition ~\ref{hg.10} with 
$M=\bbS^{k}$.
\label{hg.29}\end{proposition}
\begin{proof}
Let $T: \bbS^{k}\to \Frc(X;E)$ be a representative of an element of the
homotopy group
$\pi_{k}(\Frc(X;E))$.  We assume that $T_{s_{0}}=\Id$, where
$s_{0}\in \bbS^{k}$ is the basepoint of $\bbS^{k}$.  
By deforming $T$ if needed, we can assume that $T\equiv \Id$ in a open
ball $B^{k}_{0}\subset \bbS^{k}$ containing $s_{0}$.  Let the closed
ball $ \cb$ be the complement of $B_{0}^{k}$ in $\bbS^{k}$.

Let $V\subset H$ be a closed subspace of $H$ of finite codimension that can
be used to define the index of T:
\[
           \ind(T)= [H/V] - [H/T(V)]= [V^{\perp}]- [T(V)^{\perp}].
\]
By lemma~\ref{hg.25}, we can assume that $V^{\perp}\subset \CId(X;E)$ and
$T_{s}(V)^{\perp}\subset \CId(X;E)$ for all $s\in\bbS^{k}$.  When we restrict
$T$ to the closed ball $\cb\subset \bbS^{k}$, the index become zero because
$\cb$ is contractible and $T\equiv \Id$ on $\pa \cb$.  Thus,
by lemma~\ref{hg.25}, we can also assume that $V$ is such that 
$V^{\perp}$ and $T(V)^{\perp}$ are isomorphic vector bundles when restricted
to $\cb$.  Let $\varphi: V^{\perp}\to T(V)^{\perp}$ be a (smooth) isomorphism
between $V^{\perp}$ and $T(V)^{\perp}$ on $\cb$.  By extending $\varphi$
to act by $0$ on $V$, we get an associated family of bounded operator
$\phi: \cb\to \mathcal{L}(H,H)$.  Since $V^{\perp}\subset \CId(X;E)$ and
$T_{s}(V)^{\perp}\subset\CId(X;E)$ for all $s\in \bbS^{k}$, this is in fact
a map
\[
           \phi: \cb \to x^{\infty}\fc^{-\infty}(X;E).
\]
By construction of $\phi$ and the compactness of $\cb$, we see that for
$\lambda>0$ large enough,
\[
               T_{s}+\lambda \phi(s)\in\Gc(X;E) ,\; \forall s\in \cb.
\]
Thus, rescaling $\phi$ if needed, we can assume $T+\phi$ is a map of the form
\[
              T+\phi : \cb\to \Gc(X;E).
\]
Notice that $\no'\circ(T+\phi)=\no'\circ T$, so $T+\phi$ is in fact a lift
of the map
\[
                  \no'\circ T: \cb\to \Gso(X;E)[[x]].
\]
Moreover, since $T\equiv \Id$ on $\pa\cb\cong \bbS^{k-1}$, the map
$T+\phi$ takes value in $\Gd(X;E)$ when restricted to $\pa \cb$.  By definition
of the boundary homomorphism (see for instance section 17.1 in 
\cite{Steenrod}),
\[
    \pa( [\no'\circ T])=  [\left.(T+\phi)\right|_{\pa\cb}]\in
\pi_{k-1}(\Gd(X;E)).
\]
Now, since $V^{\perp}=T(V)^{\perp}$ canonically
on $\pa\cb$, $\phi$ is just a map
\[
             \phi: \pa\cb\to \End(V^{\perp},V^{\perp}).
\] 
Clearly, the clutching construction applied to $(\Id+\phi)^{-1}$ gives
the virtual bundle $[T(V)^{\perp}]-[V^{\perp}]$, which means that 
the clutching construction applied to $(\Id+\phi)$ gives 
$[V^{\perp}]-[T(V)^{\perp}]$.  This shows that $c\pa\nos(T)=\ind(T)$,
which concludes the proof.
\end{proof}
Since \eqref{hg.27} and \eqref{hg.28} are isomorphisms, the surjectivity
of the boundary homomorphism is equivalent to the surjectivity of the index
map of definition~\ref{hg.10}.  In \cite{Atiyah}, there is a proof of the 
surjectivity of the index map, but considering the space $\mathcal{F}(H)$
of all Fredholm operators acting on $H$.  It is possible to adapt this proof
to our situation so that it still works.  We first need to discuss further
the topology of the space $\Frc(X;E)$.

For $n\in \bbN$, consider the vector bundle $E^{n}=E\otimes\bbC^{n}$.
The bundle $E$ can be seen as a subbundle of $E^{n}$ via the inclusion
\[
    \begin{array}{rrcl}
              i: & E & \to & E\otimes \bbC^{n}  \\
                 & e  & \mapsto & 
                                   e\otimes (1,
\underbrace{0,\ldots,0}_{n-1\,\mbox{times}}) .
\end{array}          
\]  
This gives an inclusion 
\begin{equation}
\begin{array}{rrcl}
              i': & \Frc(X;E)& \hookrightarrow & \Frc(X;E^{n}) \\
                  &  \Id +A  &\mapsto & i'(\Id+A), 
\end{array}   
\label{hg.4}\end{equation}
where $i'(\Id+A)$ acts on sections of $E\subset E^{n}$ as $(\Id+A)$ do 
and acts on sections of of the complement of $E$ in
$E^{n}$ as the identity.  
\begin{definition}
A subspace $\mathcal{U}$ of a topological space $\mathcal{T}$ is said 
to be a \textbf{weak deformation retract} of $\mathcal{T}$ if for any closed
manifold $M$ and any continuous map $f:M\to \mathcal{T}$, there exists 
a continuous map $g$ homotopic to $f$ such that $g(M)\subset \mathcal{U}$.
\label{hg.5}\end{definition}
\begin{lemma}
Under the inclusion \eqref{hg.4}, $\Frc(X;E)$ is a weak deformation retract
of $\Frc(X;E^{n})$.
\label{hg.6}\end{lemma}
\begin{proof}
Let $M$ be a closed manifold and let $f:M\to \Frc(X;E^{n})$ be a continuous 
map.  Without loss of generality, we can assume $f$ is smooth.  Compose
it with the normal operator $\no$ to get a map
\[
       \widetilde{f}=\no\circ f: M\to \Gs(\pa X;E^{n}).
\]
Let $\Delta_{\pa X/Y}^{E}$ be a family of Laplacians associated to the 
fibration $\Phi:\pa X\to Y$ and the complex vector bundle $E$ (not $E^{n}$).
By the compactness of $M$, applying an argument similar to the proof of 
lemma~\ref{ckc.8}, we see that there exists a spectral section 
$(\Pi,R_{1},R_{2})$ for $\Delta_{\pa X/Y}^{E}$ such that 
$(\Pi^{n}=\bigoplus_{i=1}^{n}\Pi,R_{1},n R_{2})$, which is a spectral section
for $\bigoplus_{i=1}^{n}\Delta_{\pa X/Y}^{E}$, satisfies the estimate
\[
  \| \widetilde{f}_{b}(m) -\Id -\Pi^{n}_{b}(\widetilde{f}_{b}(m)-\Id)
\Pi^{n}_{b}\|
      < \frac{1}{2\| \widetilde{f}_{b}(m)^{-1}\|},\quad \forall b\in \fcn^{*}Y,
            m\in M.
\]
By the proof of corollary~\ref{ckc.10}, this means that
\[
  t\mapsto \widetilde{f}_{t}= \widetilde{f} + t(\Pi^{n}(\widetilde{f})\Pi^{n}-
       (\widetilde{f}-\Id)), \quad t\in[0,1],
\]
is a homotopy of smooth maps from $M$ to $\Gs(\pa X;E)$.

Let $F$ be the range of $\Pi$.  Taking $\Pi$ to have range with a larger rank
if necessary, we can assume that $\widetilde{f}_{1}$ is homotopic, through 
smooth maps from $M$ to $\Gamma_{\Id}(\fcn^{*}Y; \pi^{*}\GL(F^{n},\bbC))$,
to a map $\widetilde{f}_{2}$ with
\[
    \widetilde{f}_{2}(M)\subset \Gamma_{\Id}(\fcn^{*}Y;\pi^{*}\GL(F,\bbC)),
\]
where $\GL(F,\bbC)\subset \GL(F^{n},\bbC)$ acts as the identity on the 
complement of $F$ in $F^{n}$.  In particular, this means that
\[
 \widetilde{f}_{2}(M)\subset\Gs(\pa X;E)\subset \Gs(\pa X; E^{n}),
\]
where the inclusion $\Gs(\pa X;E)\subset \Gs(\pa X; E^{n})$ is induced
from the inclusion \eqref{hg.4}.  Let $f_{2}: M\to \Frc(X;E)$ be a map such 
that
$\widehat{\no(f_{2})}=\widetilde{f}_{2}$, which is possible by Seeley extension
for manifold with corners.  Then, using again Seeley extension, one can 
lift the homotopy between $\widetilde{f}=\widetilde{f}_{0}$ and 
$\widetilde{f}_{2}$ to a homotopy between $f$ and $f_{2}$ through
smooth maps from $M$ to $\Frc(X;E^{n})$.

\end{proof}

\begin{proposition}
For $k\in \bbN$, the boundary homomorphism
\[
          \pa : \pi_{k}(\Gso(X;E)[[x]])\to \pi_{k-1}(\Gd(X;E))
\]
is surjective.
\label{hg.30}\end{proposition}
\begin{proof}
By proposition~\ref{hg.29}, it suffices to show that the index map
\[
              \ind: \pi_{k}(\Frc(X;E))\to \widetilde{K}^{0}(\bbS^{k})
\]
is surjective.

Let $\E\to \bbS^{k}$ be an arbitrary complex vector bundle over $\bbS^{k}$.
Let $\F$ be another complex vector bundle such that 
$\E\oplus\F= \underline{\bbC}^{n}$ is trivial.  For $s\in \bbS^{k}$, let
$p_{s}\in\End(\bbC^{n})$ be the projection onto $\E_{s}\subset\bbC^{n}$.

Let $T_{-1}\in\Frc(X;E)$ be an operator such that $\ind(T_{-1})=-1$.  Such
an operator exists by corollary~\ref{im.11}.  Consider the family of operators
\[
    \widetilde{T}_{s}= T_{-1}\otimes p_{s} + \Id\otimes (1-p_{s})\in
  \Frc(X;E\otimes \underline{\bbC}^{n}), \; s\in \bbS^{k},
\] 
acting on $H\otimes\bbC^{n}$.  Clearly, $\ind(\widetilde{T})=-[\E]$.  By
lemma~\ref{hg.6}, we can deform the family $\widetilde{T}$ to a family
\[
     R:\bbS^{k}\to \Frc(X;E)\subset \Frc(X;E\otimes \underline{\bbC}^{n}).
\]
By homotopy invariance of the index, $\ind(R)=-[\E]$.  If $\rank(\E)=j$,
let $T_{j}\in \Frc(X;E)$ be an operator of index $j$, then
\[
      \ind(T_{j}\circ R)=\ind(T_{j})+\ind(R)= [j]-[\E].
\]
By construction, $T_{j}\circ R_{s}$ has index zero for all $s\in\bbS^{k}$.  
In
particular, deforming $T_{j}\circ R$ if necessary, we can assume that 
$T_{j}\circ R_{s_{0}}$ is invertible at the basepoint $s_{0}\in \bbS^{k}$.
Finally, consider the family
\[
     \widetilde{R}= (T_{j}\circ R_{s_{0}})^{-1}(T_{j}\circ R): \bbS^{k}\to 
    \Frc(X;E).
\]
By construction, $\widetilde{R}_{s_{0}}=\Id$, and 
$\ind(\widetilde{R})=[j]-[\E]$.  Since any element of 
$\widetilde{K}^{0}(\bbS^{k})$ is of the form $[j]-[\E]$, $j=\rank{\E}$,
this shows that the index map is surjective.
\end{proof}

\begin{theorem}
For $k\in \bbN_{0}$ and provided $\dim Z>0$,
\[
\begin{array}{lcl}
 \pi_{2k+1}(\Gc(X;E)) &\cong & \widetilde{K}^{-1}(Y^{T^{*}Y}), \\
  \pi_{2k}(\Gc(X;E)) &\cong & \ker\left[\ind_{t}:\K(T^{*}Y)\to \bbZ\right].
\end{array}
\]
\label{hg.31}\end{theorem}
\begin{proof}
The case of the set of connected components was  handled in 
lemma~\ref{hg.23}.  For the remaining cases, the surjectivity of 
the boundary homomorphism allows us to decompose the long exact 
sequence of homotopy groups \eqref{hg.21} into short exact sequences.
More precisely, for $k\in\bbN$, we get the exact sequences
\begin{equation}
    0 \longrightarrow \pi_{2k-1}(\Gc(X;E))\longrightarrow 
\pi_{2k-1}(\Gso(X;E))\longrightarrow 0
\label{hg.32}\end{equation}
and
\begin{equation}
0\longrightarrow \pi_{2k}(\Gc(X;E))\longrightarrow \pi_{2k}(\Gso(X;E))
\overset{\pa}{\longrightarrow} \pi_{2k-1}(\Gd(X;E))\longrightarrow 0,
\label{hg.33}\end{equation}
where we used the fact that $\pi_{2k}(\Gd(X;E))\cong\{0\}$.  For odd
homotopy groups, the theorem easily follows from \eqref{hg.32}and 
remark~\ref{hg.22}.  For even homotopy groups, notice that by 
proposition~\ref{hg.29} and proposition~\ref{hg.3}, the exact sequence
\eqref{hg.33} can be rewritten 
\begin{equation}
0\longrightarrow \pi_{2k}(\Gc(X;E))\longrightarrow \pi_{2k}(\Frc(X;E))
\overset{\ind}{\longrightarrow} \widetilde{K}^{0}(\bbS^{2k})\longrightarrow 0,
\label{hg.34}\end{equation}
so $\pi_{2k}(\Gc(X;E)) \cong  \ker\left[\ind:\pi_{2k}(\Frc(X;E))\to 
\widetilde{K}^{0}(\bbS^{k})\right]$.  But the homotopy group 
$\pi_{2k}(\Frc(X;E))$ is
isomorphic to $\K(T^{*}Y)$, 
which
is a finitely generated  $\bbZ$-module, and 
$\widetilde{K}^{0}(\bbS^{k})\cong \bbZ$.  By the classification of finitely
generated $\bbZ$-modules\footnote{See for instance in \cite{Hartley-Hawkes}.}, 
the isomorphism class of the kernel is the same
for all surjective homomorphism  $\K(T^{*}Y)\to \bbZ$.  Since 
$\ind_{t}:\K(T^{*}Y)\to \bbZ$ is also surjective, this means
\[
 \pi_{2k}(\Gc(X;E)) \cong  \ker\left[\ind_{t}:\K(T^{*}Y)\to \bbZ\right].
\]
\end{proof}

\section{Vanishing of the Homotopy Groups}\label{vhg.0}

In the case of cusp operators, that is, when the fibration 
$\Phi:\pa X\to \{\pt\}$
is trivial, one gets as an easy corollary of theorem~\ref{hg.31} the 
weak contractibility result of \cite{fipomb}, namely that all
the homotopy groups of $\Gc(X;E)$ vanish.  

\begin{corollary}
If $\Phi:\pa X\to \{\pt\}$ is the trivial fibration (cusp operators), then
$\Gc(X;E)$ is weakly contractible, meaning that all its homotopy groups are
trivial, including the set of connected components.
\label{vhg.1}\end{corollary}
\begin{proof}
In this case, $Y=\pt$, so $T^{*}Y$ does not really make sense, but looking back
at the identification \eqref{ckc.32} and the proof of lemma~\ref{hg.2},
we see that $\K(T^{*}Y)$ and $\widetilde{K}^{-1}(Y^{T^{*}Y})$ must be replaced
by $\widetilde{K}^{0}(\bbS^{2})\cong \bbZ$ and $\widetilde{K}^{-1}(\bbS^{2})
\cong\{0\}$ in the statement of theorem~\ref{hg.31}.  The result then
easily follows from the fact that 
$\ind_{t}:\widetilde{K}^{0}(\bbS^{2})\to \bbZ$ 
is an isomorphism.
\end{proof}

In fact, it turns out that with 
only little more effort, one can deduce that $\Gc(X;E)$ is actually 
contractible.  It suffices to apply with slight modifications the argument of 
Kuiper in his proof of
the contractibility of the group of invertible bounded operators acting
on a separable Hilbert space (see \cite{Kuiper}).   
   
The main difference in our situation is that the topology we consider
on the group $\Gc(X;E)$ is not the one coming from the operator norm, but 
the $\CI$-topology coming from the Schwartz kernels of the smoothing 
operators in $\fc^{-\infty}(X;E)$.  Nevertheless, since
$\Gc(X;E)\subset \fc^{0}(X;E)$, we can also give to $\Gc(X;E)$ the 
toplogy induced by the operator norm $\| \cdot \|$ of bounded operators 
acting on $\Ld(X;E)$.  This is a weaker topology then the $\CI$-infinity
topology, in the sense that it has less open sets.  In the operator norm 
topology, the space $\Gc(X;E)$ is easily seen to be a metric space.  By
a theorem of Stone\footnote{See for instance theorem $4$, section
$I.8.4$, p.101 in \cite{schubert}}, 
metric spaces are paracompact.  This will be very
useful to retract  $\Gc(X;E)$ (in the $\CI$-topology) onto a CW-complex.
   
\begin{definition}
We will say that an open ball in $\Id+\fc^{-\infty}(X;E)$ of radius $\epsilon$
\[
     B_{\epsilon}(A)=\{ \Id + Q \; | \;Q\in \fc^{-\infty}(X;E),\;\; 
                          \|A-\Id-Q\|<\epsilon \}, \quad A\in \Gc(X;E), 
\]
is \textbf{small} if $B_{3\epsilon}(A)\subset \Gc(X;E)$.  Clearly, one can
cover $\Gc(X;E)$ by such balls.
\label{vhg.2}\end{definition}

\begin{proposition}
If $\Gc(X;E)$ is weakly contractible, then it is contractible.  In particular,
if $\Phi:\pa X\to \pt$ is the trivial fibration, then $\Gc(X;E)$ is 
contractible.
\label{vhg.3}\end{proposition}
\begin{proof}
Let $\{B_{\epsilon_{i}}(A_{i})\}_{i\in I}$ be a covering of $\Gc(X;E)$ by
small balls,
\begin{equation}
              \Gc(X;E)=\bigcup_{i\in I} B_{\epsilon_{i}}(A_{i}).
\label{vhg.4}\end{equation}
Since the Banach space of bounded operators acting on $\Ld(X;E)$ is 
separable, we can assume $I=\bbN$.  Moreover, since $\Gc(X;E)$ is 
paracompact in the operator norm topology, we can assume that the covering
\eqref{vhg.4} is locally finite and that it has an associated partition
of unity $\{\phi_{i}\}_{i\in\bbN}$.  A priori, the $\phi_{i}$ are continuous
with respect to the operator norm topology, but since the $\CI$-topology is 
a finer topology, this means they are also continuous with respect to the 
$\CI$-topology.  For $t\in [0,1]$, consider the
following homotopy
\begin{equation}
\begin{array}{rrcl}
       \xi_{t}:& \Gc(X;E) & \to & \Gc(X;E) \\
               & z & \mapsto & (1-t)z + t\underset{i\in\bbN}{\sum}
         \phi_{i}(z)A_{i}.
\end{array}
\label{vhg.6}\end{equation}
For $t=0$, $\xi_{0}$ is just the identity.  To see that 
the image really lies in $\Gc(X;E)$, let $z\in\Gc(X;E)$ be given.  Then
there exists a neighborhood $\mathcal{U}$ of $z$ such that $\mathcal{U}$ has
a non-empty intersection with only finitely many open sets of the covering
\eqref{vhg.4}, say $B_{\epsilon_{i_{1}}}(A_{i_{1}}),\ldots,
B_{\epsilon_{i_{m}}}(A_{i_{m}})$.  Without loss of generality,
assume that $B_{\epsilon_{i_{1}}}(A_{i_{1}}),\ldots,
B_{\epsilon_{i_{m}}}(A_{i_{n}})$, $n\le m$ are the open sets of the covering
containing $z$, and assume that 
$\epsilon_{i_{1}}=\max\{\epsilon_{i_{1}},\ldots, \epsilon_{i_{n}}\}$.  By 
construction, we have
\begin{equation}
  z\in B_{\epsilon_{i_{k}}}(A_{i_{k}})\subset B_{3\epsilon_{i_{1}}}(A_{i_{1}}),
        \quad \forall k\in\{1,\ldots,n\},  
\label{vhg.5}\end{equation}
which implies that $\xi_{t}(z)\in B_{3\epsilon_{i_{1}}}(A_{i_{1}})\subset 
\Gc(X;E)$ for all $t\in [0,1]$.  

Let $N$ be the nerve of the covering \eqref{vhg.4}.  Call $b_{i}$ the vertex 
that corresponds to the open set $B_{\epsilon_{i}}(A_{i})$.  Then $N$ is a
CW-complex with affine simplices as cells.  There is an obvious inclusion
$\rho:N\to \Gc(X;E)$ given by sending the vertex $b_{i}$ to $A_{i}$ for
$i\in \bbN$ and so that on any simplex of $N$, $\rho$ is an affine
map.  Since $\xi_{1}(\Gc(X;E))\subset N$, the homotopy \eqref{vhg.6} shows 
that the inclusion $\rho$ is a homotopy equivalence.  By Whitehead's theorem,
we then conclude that $\Gc(X;E)$ is contractible.

\end{proof}

\providecommand{\bysame}{\leavevmode\hbox to3em{\hrulefill}\thinspace}
\providecommand{\MR}{\relax\ifhmode\unskip\space\fi MR }
\providecommand{\MRhref}[2]{%
  \href{http://www.ams.org/mathscinet-getitem?mr=#1}{#2}
}
\providecommand{\href}[2]{#2}

\end{document}